\documentclass[]{interact}
\usepackage{epstopdf}
\usepackage[caption=false]{subfig}
\usepackage{amsmath,amsfonts,amssymb,adjustbox}
\usepackage[sort&compress,square,comma,numbers]{natbib}
\bibpunct[, ]{[}{]}{,}{n}{,}{,}
\usepackage{xcolor,hyperref}
\usepackage{enumitem,kantlipsum}
\usepackage[normalem]{ulem}
\newcommand{\hluo}[1]{{\color{black}#1}}
\newcommand{\hluoold}[1]{{\color{blue}\sout{#1}}}
\newcommand{\steve}[1]{{\color{black}#1}}

\newcommand{\mario}[1]{{\color{black}#1}}
\newcommand{\hrluonew}[1]{{\color{black}#1}}

\makeatletter
\def\NAT@def@citea{\def\@citea{\NAT@separator}}
\makeatother

\theoremstyle{plain}
\newtheorem{theorem}{Theorem}[section]
\newtheorem{lemma}[theorem]{Lemma}
\newtheorem{corollary}[theorem]{Corollary}
\newtheorem{proposition}[theorem]{Proposition}

\theoremstyle{definition}
\newtheorem{definition}[theorem]{Definition}
\newtheorem{example}[theorem]{Example}

\theoremstyle{remark}

\usepackage{amsthm}
\usepackage{mathtools}

\usepackage{lmodern}
\usepackage{color}
\usepackage[normalem]{ulem}
\usepackage{placeins}
\usepackage{hyperref}
\usepackage{graphicx}
\usepackage{adjustbox} 

\providecommand{\tabularnewline}{\\}
\definecolor{burgundy}{rgb}{0.5, 0.0, 0.13}

\makeatother

\usepackage[english]{babel}
\usepackage{setspace}
\begin{document}
\articletype{}

\title{Asymptotics of Lower Dimensional Zero-Density Regions}

\author{
\name{Hengrui Luo\textsuperscript{a}\thanks{CONTACT Hengrui Luo. Email: hrluo@lbl.gov}, Steven N. MacEachern\textsuperscript{b} and Mario Peruggia\textsuperscript{b}}
\affil{\textsuperscript{a}Department of Statistics, College of Arts and Sciences, the Ohio State University, 1958 Neil Ave. Columbus OH 43210 USA; and Lawrence Berkeley National Laboratory, Berkeley, CA 94720, USA; \textsuperscript{b}Department of Statistics, College of Arts and Sciences, the Ohio State University, 1958 Neil Ave. Columbus OH 43210, USA}
}

\maketitle

\begin{abstract}
Topological data analysis (TDA) allows us to explore the topological features of a dataset. Among topological features, lower dimensional
ones have recently drawn the attention of practitioners in mathematics and statistics due to
their potential to aid the discovery of low dimensional structure
in a data set. However, lower dimensional features are usually challenging
to detect 
\mario{based on finite samples and using TDA methods that ignore the probabilistic mechanism that generates the data.}
In this paper, lower dimensional topological features occurring
as zero-density regions of density functions are introduced and thoroughly
investigated. Specifically, we 
consider sequences of coverings for the support of a density function in which the coverings are comprised of balls with shrinking radii.
We show that, when these coverings satisfy certain sufficient conditions as the sample size goes to infinity, we can detect lower dimensional, zero-density regions with increasingly higher probability while guarding against false detection. We supplement the theoretical developments with
the discussion of simulated experiments that elucidate the behavior of the methodology for different choices of the tuning parameters
that govern the construction of the covering sequences and characterize the asymptotic results.
\end{abstract}



\section{Introduction}
\mario{
We consider the problem of identifying certain lower-dimensional, geometric features of (the support of) a density that generates a stream of observed data. What makes this identification at once interesting and challenging is that the features we are concerned with cannot be detected using any finite amount of data. Typical topological data analysis (TDA) techniques would be unsuitable in this context because they rely on the relative arrangement of points in finite data sets to understand the structure of the (sub)-space where the data occur (see for example 
\cite{Carlsson2009,carlsson2021topological} and  \cite{Edelsbrunner&Harer2010} for introductory treatments).} 
\mario{We can make progress in the detection process because we rely on knowledge of relevant distributional properties of the data generating mechanism that we acquire as more and more data accrue. This is in contrast to a standard use of TDA methods that would not explicitly account for these distributional properties.
We study the asymptotics of the data generating mechanism directly and show that they can be linked to certain topological characteristics. 
}

\begin{figure}[t!]
\begin{center}
\includegraphics[width=8cm,height=6cm]{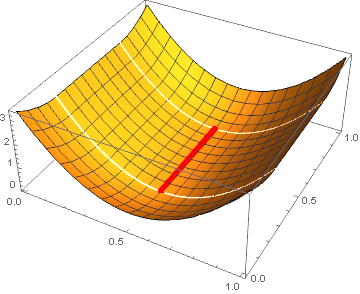}
\end{center}
\caption{\label{Minimal Example}A lower dimensional region $S_{0}=\{\frac{1}{2}\}\times[0.25,0.75]$
(red segment) for the density function $f(\boldsymbol{x})\propto d(\boldsymbol{x},S_{0})^{4}\circ\boldsymbol{1}_{[0,1]^{2}}(\boldsymbol{x})$,
where $d(\boldsymbol{x},S_{0})=\inf_{\boldsymbol{y}\in S_{0}}d(\boldsymbol{x},\boldsymbol{y})$
and $d$ denotes Euclidean distance. The density $f$ evaluates to
0 over $S_{0}$, but nowhere else over $[0,1]^{2}$.}
\end{figure}

In this paper, we first discuss situations where independent and identically distributed (i.i.d.) data points are drawn from  a distribution having \hluo{a continuous} density function $f$ (with respect to Lebesgue measure on $\mathbb{R}^{d}$)  on  its support $supp(f)=M\subset\mathbb{R}^{d}$. Our results are first
formally stated for $supp(f)=M=[0,1]^{d},$ $d\in\mathbb{N}$, and
then extended to more general situations. We work with a well-behaved
version of the density $f$ for which the notion of a zero-density
region $S_{0}\subset supp(f)\subset M$ (to be formally defined later) is meaningful.
Such a zero-density region $S_{0}$ is difficult to identify with
traditional constructions of simplices or density estimators. 
\begin{example}
Let $S_{0}=\{\frac{1}{2}\}\times[0.25,0.75]$ and define $d(\boldsymbol{x},S_{0})=\inf_{\boldsymbol{y}\in S_{0}}d(\boldsymbol{x},\boldsymbol{y})$,
where $d$ denotes the $L^2$ metric. Consider the density $f(\boldsymbol{x})\propto d(\boldsymbol{x},S_{0})^{4}\circ\boldsymbol{1}_{[0,1]^{2}}$
shown in Figure \ref{Minimal Example}, for which $S_{0}$ is a zero
density region of lower dimension. $S_{0}$ does not contain any probability
mass and, being a segment, is a ``lower dimensional object'' (a
concept to be made more precise later). 

The volume of $S_{0}$ as a subset of the support of the density, $[0,1]^2$, is zero and the density assigns positive probability to \steve{every} 
neighborhood of \steve{every} 
point in $S_{0}$. \steve{Sampled points ``on either side'' of $S_0$ can be arbitrarily close to one another.}  Hence,   
it 
\steve{might seem to} be impossible to identify the topological structure of $S_{0}$ with accuracy.  \steve{However, consideration of an asymptotic argument that relies on }
the rate at which points accumulate in 
\steve{the vicinity of $S_{0}$} as the sample size grows \steve{allows us to identify the structure of $S_{0}$}.
\steve{We note that the} situation would be different if the density were zero on a region of positive volume, such as a disk, and nonzero elsewhere.  The disk could then be identified with a large enough sample as a ``hole'' not containing any points.

Our upcoming results 
\steve{consider samples of increasing size.  The main result}
allows us to detect the lower dimensional object $S_{0}$ with probability one as the sample size $n$
goes to infinity\steve{, while at the same time helping us }
avoid detection of false holes. 
\end{example}

Our approach to detecting zero-density regions is to conduct an analysis
as we vary the radius of a collection of covering balls of $supp(f)$.
By choosing the shared radius of the covering balls appropriately
relative to increasing sample size $n$, we wish that $S_{0}$
be covered only by balls having no observations inside. For each point
in the non-zero density region we wish for the point to eventually
be covered by a ball with at least one observation inside. If our
wishes come true, then we can simply collect the empty covering balls
and recover an approximation to the region $S_{0}$. In Theorem~\ref{Main Theorem}
we present a set of sufficient conditions for  our detection method to work asymptotically.

This notion of varying the radius of the covering balls can be related
to the construction of complexes in TDA. For the \v{C}ech complex,  balls
are centered at the observed points. Balls of a fixed radius $r$
lead to a \v{C}ech complex $C(\mathcal{X},r)$. Varying the radius $r$, the \v{C}ech
filtration, a collection of \v{C}ech complexes, is produced. For a small
radius $r$, the balls will not overlap and no holes will be discovered.
For a large radius $r$, lower dimensional zero-density regions will
be covered and these holes will not be found.

The rest of the paper is organized as follows. We first illustrate
our observations above with a simple example in Section~\ref{sec:An-Illustrative-Example}.
Then, in Section~\ref{sec:Statement-of-the-main-result}, we discuss our approach for the construction of coverings of
a compact support along with a set of sufficient conditions that ensure
asymptotic consistency of the procedure for the detection of zero-density regions.
Generalizations of the results to the case of a non-compact support are provided
in Section~\ref{sec:Non-compact-support}.
We present some experimental results and connections to other areas in Section~\ref{sec:Simulations-and-Connections} and end with a brief discussion of our findings in Section~\ref{sec:Conclusion}.

\section{\label{sec:An-Illustrative-Example}An Illustrative Example}

The central problem we investigate in this paper is the detection of a lower dimensional zero-density region. Before formally addressing the problem, we present an illustrative example to show how dimensionality plays an important role.

We consider three different densities on the
interval $M=[-1,1]$, $\dim M=1$, and the problem of detecting interesting topological features
by partitioning $M$ into the union of disjoint, equally sized bins. \steve{The feature (or lack thereof) for two of the densities is easily found.}  The density 
$g(x)=\frac{2}{3}\left(\boldsymbol{1}_{[-1,-1/4]}(x)+\boldsymbol{1}_{[1/4,1]}(x)\right)$
has a genuine hole, $(-1/4,1/4)$, consisting of an interval of dimension one. This hole is
apparent, as there will never be any observations in it\steve{, nor in bins contained in it}. As the sample
size $n$ goes to infinity, one need only consider bins whose widths
decrease at a rate no faster than $n^{-1+\varepsilon}$, $\varepsilon>0$,
to ensure that the bins away from the hole are eventually filled.
When $S_{0}=\emptyset$, as is the case for the density $h(x)=\frac {3}{8}\boldsymbol{1}_{[-1,1]}(x)\cdot\left(x^{2}+1\right)$,
and the width of the bins decreases at a rate no faster than $n^{-1+\varepsilon}$, \steve{all of the bins will eventually be filled and it will be evident that there is no topological feature.}  

The interesting case has $S_{0}\neq\emptyset$ and $\dim S_{0}<\dim M=1$,
as is the case for $f(x)=\frac{3}{2}x^{2}\boldsymbol{1}_{[-1,1]}(x)$. Here,
$S_{0}=\{0\}$, and the question is whether we can detect this topological
feature of dimension zero with positive probability. The hole is difficult
to detect because samples will accumulate \steve{in every neighborhood of}
~$S_{0}=\{0\}$. If
the bins \steve{shrink too slowly,} 
\steve{the sequence of bins that contain $\{0\}$ will be filled} 
as \hluo{the sample size} $n$ goes to infinity.  

\begin{figure}[t]
\begin{adjustbox}{center}
    \begin{tabular}{ccc}
    $f(x)=\frac{3}{2}x^{2}\boldsymbol{1}_{[-1,1]}(x)$ & $g(x)=\frac{2}{3}\left(\boldsymbol{1}_{[-1,-\frac{1}{4}]}(x)+\boldsymbol{1}_{[\frac{1}{4},1]}(x)\right)$ & $h(x)=\frac{3}{8}\boldsymbol{1}_{[-1,1]}(x)\cdot\left(x^{2}+1\right)$\tabularnewline
    \includegraphics[scale=0.25]{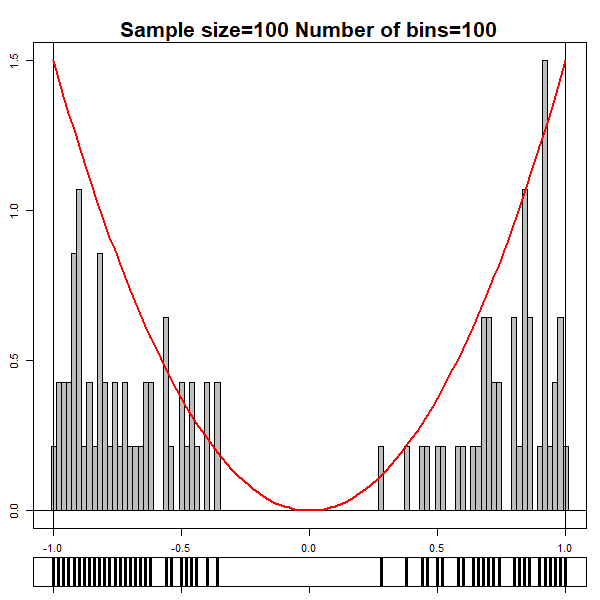} & \includegraphics[scale=0.25]{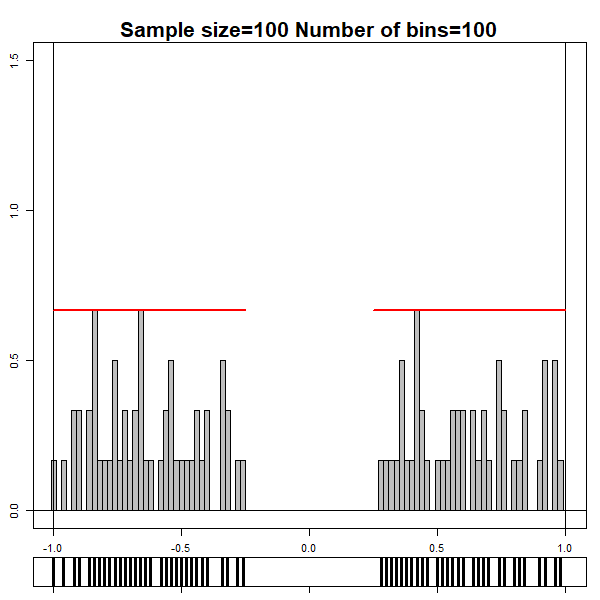} & \includegraphics[scale=0.25]{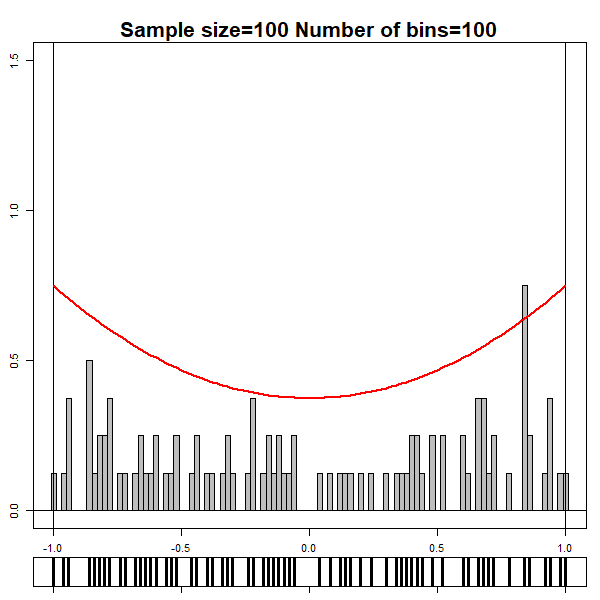}\tabularnewline
    \includegraphics[scale=0.25]{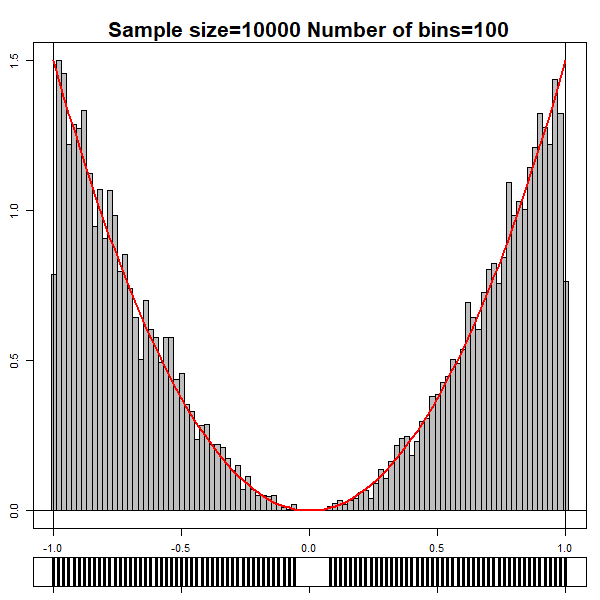}  & \includegraphics[scale=0.25]{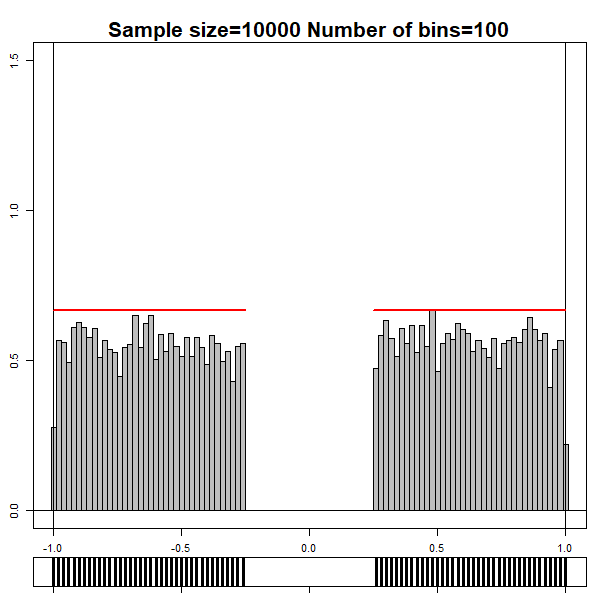} & \includegraphics[scale=0.25]{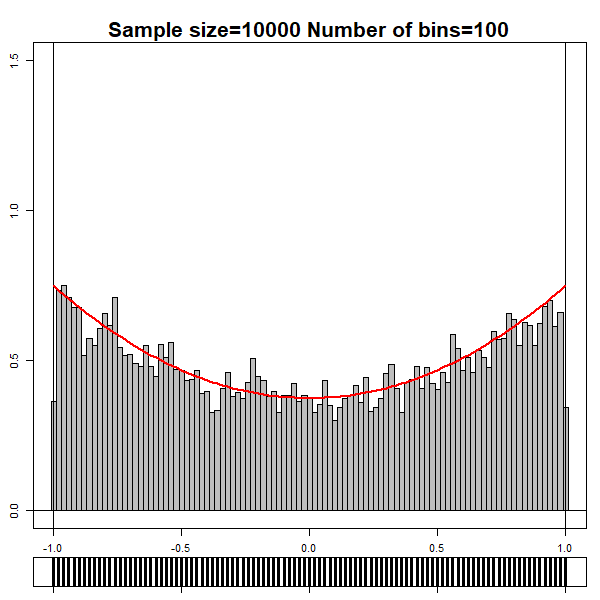} \tabularnewline
    \end{tabular}
\end{adjustbox}
\caption{\label{fig:The-binning-scheme} Binned histogram and binary heatmaps
based on random samples of size $n=100$ and $n=10,000$ from three
distributions with densities $f$, $g$, and $h$.  All graphical summaries are constructed using
100 equal-width bins over the interval $[-1,1]\subset\mathbb{R}$.\protect \\
}
\end{figure}

In 
Figure \ref{fig:The-binning-scheme}, the binary heatmaps show
the denseness of the non-empty bins. The first column represents histograms and heatmaps for the density $f$
with a zero-dimensional hole $\{0\}$. The second column represents histograms and heatmaps for the density $g$ with a one-dimensional hole $(-1/4,1/4)$. The third column represents histograms and heatmaps for the density $h$
with no holes. \steve{The top row is for a small sample size; the bottom row for a larger sample size.}

\steve{The figure shows the ease with which binning identifies densities bounded away from $0$.  See the filled bins under the red lines in the second and third columns.  The figure also shows the ease with which a hole of full dimension is found.  See the empty bins in the center of the plots in column two.  The difficult case appears in the first column.  Here, we see a lower-dimensional hole in the density.  }
The figure suggests that, as the sample size $n$ goes to infinity, we should be able to 
\steve{detect this hole} by using binary heatmaps with an appropriate scaling scheme.
This is perhaps surprising, because \steve{the density} $f$ 
\steve{in the first column has full} support.  Formally,
our Theorem \ref{Main Theorem} shows that, if we shrink the common
width of the bins (which corresponds to the radius of the covering
balls in higher dimensions) at an appropriate rate as a function of
sample size $n$, the lower dimensional 
holes will be characterized in terms
of empty bins with probability tending to one. This result can be
extended to more general situations.

The objective of this article is to show that, without appropriate
scaling, authentic topological holes of strictly lower dimension
cannot be detected, while, with appropriate scaling, they can. 
The sufficient conditions that we will impose on the scaling schemes to attain these
results depend on the dimension of the zero-density region $S_{0}$
and also on the local smoothness of the density that generates the
data.

\section{\label{sec:Statement-of-the-main-result}\steve{Main Result for Compact Support}}

\noindent Consider a random vector $X$ having density~$f$ with respect
to Lebesgue measure on $\mathbb{R}^{d}$. For any Borel set $A\in\mathbb{R}^{d}$,
let $P(X\in A)=\int_{A}f(x)\,dx$. 
\begin{definition}
\label{def:(Support)-We-assume}(Support) The \emph{support} of $X$
is defined to be 
\begin{align}
supp(X) & \coloneqq\cap_{\{R_{X}\text{ closed in }\mathbb{R}^{d}\mid P(X\in R_{X})=1\}}R_{X}.
\end{align}
With an abuse of notation, we write $supp(f)$ to represent the support
of $X$. 
\end{definition}

In most applications where TDA is employed, a compact support $M$
can be reasonably assumed. For technical convenience, we will also
assume the existence of a continuous version of the density $f$ on its
support $M$. The latter is a mild condition satisfied for most theoretical
questions and applied scenarios. In the following discussion, as will
be stated in \hyperlink{Assumption1}{Assumption~1} on page~\pageref{PageAss1-5},
we consider the case of $M=[0,1]^{d}\subset\mathbb{R}^{d},d>0$. With
this particular $M$ in mind, we turn to a special region called the
zero-density region \hrluonew{assumed to lie} in its interior. 
\begin{definition}
\label{def:(Zero-density-region)}(Zero-density region) For the continuous
version of a density $f$ on $M$, we call the inverse image of 
$\{0\}$, i.e.
$f^{-1}(\{0\})$, the \emph{zero-density region} of $f$ and denote
it by $S_{0}$. 
\end{definition}

\hrluonew{
In a first result (Theorem~\ref{Main Theorem}), we suppose
that $S_0$ consists of a single, connected component (as in \hyperlink{Assumption2}{Assumption 2}) for simplicity. 
 We then extend the initial result to the case where $S_{0}$ consists of a finite number
of connected components in Corollary~\ref{cor:Suppose-we-have}.
}

The local behavior of the density $f$ around the zero-density region
$S_{0}\subset M$ is a crucial aspect of our investigation. The following
concepts will help us to describe the behavior of~$f$ near $S_{0}$. 
We consider the $L^2$ metric and the associated norm, $\|\cdot\|$, in the following discussion, but remark that our methods can be generalized to other metrics and norms.

\noindent A \emph{ball} of radius $r$ centered at $\boldsymbol{x}$
in $\mathbb{R}^{d}$ is the open set 
\begin{align*}
B_{r}(\boldsymbol{x})\coloneqq\left\{ \left.\boldsymbol{y}\in\mathbb{R}^{d}\,\right|\|\boldsymbol{x}-\boldsymbol{y}\|<r\right\}  & .
\end{align*}
An $\epsilon$-neighborhood ($\epsilon>0$) of $S_{0}$ is the open
set 
\begin{align*}
B_{\epsilon}(S_{0})\coloneqq\left\{ \boldsymbol{x}\in\mathbb{R}^{d}\left|\,\inf_{\boldsymbol{y}\in S_{0}}\|\boldsymbol{x}-\boldsymbol{y}\|<\epsilon\right.\right\}  & .
\end{align*}
(Note that the definition of an an $\epsilon$-neighborhood makes
sense for an arbitrary set, not just a zero-density region. In particular,
the $\epsilon$-neighborhood of a point $\boldsymbol{y}$, $B_{\epsilon}(\boldsymbol{y})$,
is just the open ball of radius $\epsilon$ centered at $\boldsymbol{y}$.)
\begin{definition}
\label{def:(Order-of-smoothness)-1}(Upper and lower $\epsilon$-order
of smoothness) We consider densities $f$ for which the following
quantities, called the \emph{upper} and \emph{lower $\epsilon$-order
of smoothness,} respectively, exist and are well-defined for every
$0<\epsilon<1$: 
\[
\overline{K_{f}}(\epsilon)\coloneqq\inf\left\{ \alpha>0\left|0<\inf_{\boldsymbol{x}\in B_{\epsilon}(S_{0})\backslash S_{0}}\frac{f(\boldsymbol{x})}{d(\boldsymbol{x},S_{0})^{\alpha}}\right.\right\} ,
\]
and 
\[
\underline{K_{f}}(\epsilon)\coloneqq\sup\left\{ \alpha>0\left|\infty>\sup_{\boldsymbol{x}\in B_{\epsilon}(S_{0})\backslash S_{0}}\frac{f(\boldsymbol{x})}{d(\boldsymbol{x},S_{0})^{\alpha}}\right.\right\} .
\]
\end{definition}

Because our result concerns the limiting case as $\epsilon$ goes
to zero, we give the following definition.
\begin{definition}
\label{def:(Order-of-smoothness)}(Upper and lower order of smoothness)
The \emph{upper }and \emph{lower order of smoothness} of the density
$f$ w.r.t.\ $S_{0}$, are 
\begin{align*}
\overline{K_{f}} & =\lim_{\epsilon\rightarrow0^{+}}\overline{K_{f}}(\epsilon)\,\,\,\text{and}\,\,\,\underline{K_{f}}=\lim_{\epsilon\rightarrow0^{+}}\underline{K_{f}}(\epsilon),
\end{align*}
respectively, provided they exist. 
And if they exist, $\overline{K_{f}}\geq \underline{K_{f}}$ by this definition.
\end{definition}

\label{PageK_eps} It can be shown that $\overline{K_{f}}$ and $\underline{K_{f}}$
will exist if $\overline{K_{f}}(\epsilon)$ and $\underline{K_{f}}(\epsilon)$
exist for some $\epsilon$ and that there are densities for which
$\overline{K_{f}}(\epsilon)$ and $\underline{K_{f}}(\epsilon)$ are
undefined for all $\epsilon$. The limits $\overline{K_{f}}$ and
$\underline{K_{f}}>0$, if they exist, need not coincide. If they
do exist and coincide, we will use the notation $K_{f}=\overline{K_{f}}=\underline{K_{f}}$. 
Similar but more specific variants of this concept of order of smoothness are studied both in TDA \cite{guibas2013witnessed} as growth dimensions and in the density estimation literature as local smoothness \cite{Mammen&Polonik2013}. The generality of this kind of description makes our \hyperlink{allAssumptions}{assumptions A.1-A.6} for the main theorem (especially A.3-A.4) relevant to research from both communities. 

\begin{figure}[t!]
\begin{center}
\includegraphics[width=8cm,height=6cm]{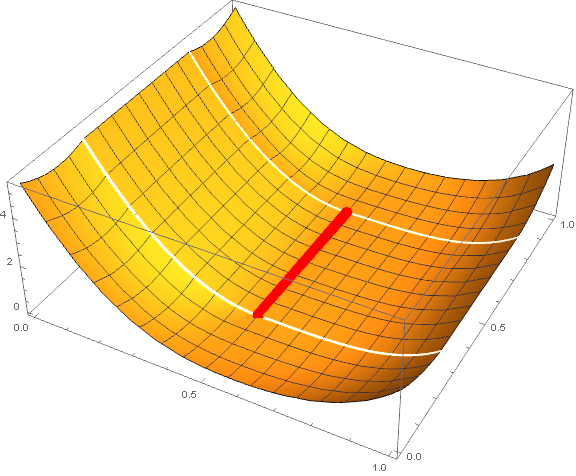}
\end{center}
\caption{\label{Minimal Example-1} Lower dimensional region $S_{0}=\{\frac{1}{2}\}\times[0.25,0.75]$
of a density function $f(\boldsymbol{x})$ with $\overline{K_{f}}\protect\neq\underline{K_{f}}$.
The density $f$ evaluates to 0 over $S_{0}$ but nowhere else over
$[0,1]^{2}$.}
\end{figure}

 The following example, depicted in Figure~\ref{Minimal Example-1},
shows that the upper and lower orders of smoothness $\overline{K_{f}}$
and $\underline{K_{f}}$ need not coincide.

\FloatBarrier
\begin{example}
\label{exa:ord of smoothness example} For $\boldsymbol{x}=(x,y)\in\mathbb{R}^{2}$,
define $h(\boldsymbol{x})$ as 
\[
h(\boldsymbol{x})=\begin{cases}
d(\boldsymbol{x},S_{0})^{4}, & \boldsymbol{x}\in[0.5,1]\times[0.25,0.75],\\
d(\boldsymbol{x},S_{0})^{2}, & \boldsymbol{x}\in[0,0.5)\times[0.25,0.75],\\
d(\boldsymbol{x},S_{0})^{4-\frac{2}{\pi}\theta_{1}}, & \boldsymbol{x}\in[0,1]\times(0.75,1],\\
d(\boldsymbol{x},S_{0})^{4+\frac{2}{\pi}\theta_{2}}, & \boldsymbol{x}\in[0,1]\times[0,0.25),
\end{cases}
\]
where $d$ is the  $L^2$ metric, $\theta_{1}=\arctan\{(y-0.75)/(x-0.5)\}$, and $\theta_{2}=\arctan\{(y-0.25)/(x-0.5)\}$.
Consider the density function $f$ satisfying the relationship $f(\boldsymbol{x})\propto g(\boldsymbol{x})\times\boldsymbol{1}_{[0,1]^{2}}(\boldsymbol{x})$.
For $\boldsymbol{x}\in M$ with $d(\boldsymbol{x},S_{0})=\delta$,
$\delta>0$, $\sup_{\boldsymbol{x}}f(\boldsymbol{x})=\delta^{2}$
and $\inf_{\boldsymbol{x}}f(\boldsymbol{x})=\delta^{4}$. Thus $\overline{K_{f}}(\epsilon)=4$
and $\underline{K_{f}}(\epsilon)=2$. 
\end{example}
\subsection{Covering balls and dimension of \texorpdfstring{$S_{0}$}{S0}}
We cover the support $M$ of the density $f$ with a collection of balls of equal radius and, letting the radius shrink to zero at an appropriate rate, we attempt to detect the zero-density region $S_{0}$
with a certain limiting probability guarantee. In what follows we will make use of the following piece of notation. 
\newline \noindent \textbf{Notation.} (Covering) Let $E$ be a given subset of $\mathbb{R}^{d}$.
We denote by $\mathcal{B}_{r}^{d}(E)$ a collection of $d$-dimensional
balls of radius $r$ whose union contains $E$.  Note that a covering $\mathcal{B}_{r}^{d}(E)$ may also depend on
the sample size $n$ through the radius $r=r(n)$ in subsequent developments. We also use $|\mathcal{B}_{r}^{d}(E)|$ to denote the cardinality of $\mathcal{B}_{r}^{d}(E)$ and write $\mathcal{B}_{r}^{d}(E)=\mathcal{B}$ when the meaning is unambiguous.

We distinguish between three types of covering balls based on how far they are from the zero-density region $S_0$. 
\begin{definition}
\label{def:(Different-types-of-covering-balls)}(Types of covering
balls) Consider a density $f$ with respect to Lebesgue measure on
$\mathbb{R}^{d}$ that is continuous on $M$ and zero on $\mathbb{R}^{d}\backslash M$.
Let $S_{0}$ be a zero-density region for $f$. Denote by $\mathcal{B}$
a covering of $M$ with balls of radius $r$ and let $B_{r}(\boldsymbol{x})\in\mathcal{B}$.
We classify $B_{r}(\boldsymbol{x})$ into one of these three types: 
\end{definition}

\begin{figure}[t!]
\centering \includegraphics[scale=0.5]{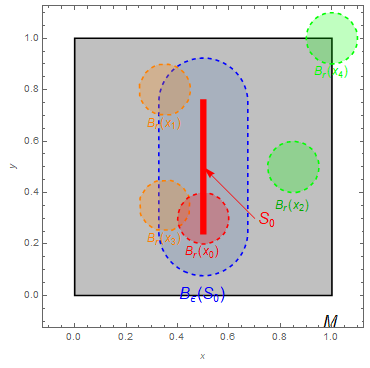}

\caption{\label{Different kinds of balls} Illustration of the types of balls
with the density in Figure~\ref{Minimal Example}. 
The blue region is an epsilon neighborhood surrounding the zero-density
region $S_{0}$. 
The red ball is an $\epsilon$-inside ball, the orange balls are $\epsilon$-neighboring
balls and the green balls are $\epsilon$-outside balls.}
\end{figure}

\begin{enumerate}
\item an \emph{$\epsilon$-outside ball}: a ball $B_{r}(\boldsymbol{x})$
such that $\boldsymbol{x}\notin B_{\epsilon}(S_{0})$; 
\item an \emph{$\epsilon$-neighboring ball:} a ball $B_{r}(\boldsymbol{x})$
such that $\boldsymbol{x}\in B_{\epsilon}(S_{0})$ and $B_{r}(\boldsymbol{x})\cap S_{0}=\emptyset$; 
\item an \emph{$\epsilon$-inside ball}: a ball $B_{r}(\boldsymbol{x})$
such that $\boldsymbol{x}\in B_{\epsilon}(S_{0})$ and $B_{r}(\boldsymbol{x})\cap S_{0}\neq\emptyset$. 
\end{enumerate}
The main result will rely on the condition $r\leq\epsilon/2$ to ensure
that every $\epsilon$-inside ball is contained in $B_{\epsilon}(S_{0})$.
The various types of covering balls are illustrated in Figure~\ref{Different kinds of balls}.

\begin{definition}
(Big O notation \cite{Burgisser&Cucker2013}) Given two positive
real sequences $f(n)$ and $g(n)$, we write $f\sim O(g)$ and say
that $f$ is big $O$ of $g$ if there exist constants $L_{1}$ and
$L_{2}\in(0,\infty)$ and an $n_{0}\in\mathbb{N}$, such that 
\begin{align*}
{\displaystyle L_{1}\cdot g(n)\leq f(n)\leq L_{2}\cdot g(n)} & ,\,\text{for all}\,\,n>n_{0}.
\end{align*}
\end{definition}

\noindent This condition means that the asymptotic behaviors of $f$
and $g$ are comparable.

Next, we introduce the notion of Minkowski dimension which we will
later use to characterize the dimension of $S_{0}$. 
\begin{definition}
\label{def:(Minkowski-dimension,-or}(Minkowski dimension, or box-counting
dimension. Definition 3.1 in \citealp{Falconer2004}) The\emph{ upper
and lower Minkowski dimensions} of a bounded subset $E\subset M$
of $\mathbb{R}^{d}$, are defined respectively as 
\begin{align*}
\overline{\dim}_{M}(E) & \coloneqq\lim\sup_{\Delta\rightarrow0}\frac{\log\mathcal{N}_{\Delta}(E)}{-\log\Delta},\\
\underline{\dim}_{M}(E) & \coloneqq\lim\inf_{\Delta\rightarrow0}\frac{\log\mathcal{N}_{\Delta}(E)}{-\log\Delta},
\end{align*}
where $\mathcal{N}_{\Delta}(E)$ is the \emph{$\Delta$-covering number}
of $E\subset\mathbb{R}^{d}$ 
\begin{align*}
\mathcal{N}_{\Delta}(E) & \coloneqq\min\left\{ k\in\mathbb{N}\mid E\subset\cup_{i=1}^{k}B_{\Delta}(x_{i})\text{ for some }x_{i}\in\mathbb{R}^{d}\right\} ,
\end{align*}
i.e., the smallest number of balls of radius $\Delta>0$ needed to
cover $E$. We call this covering of minimal cardinality a \emph{minimal
$\Delta$-covering} for $E$. When $\overline{\dim}_{M}(E)=\underline{\dim}_{M}(E)=d_{M}$,
we define the \emph{Minkowski dimension }(or box-counting dimension)
$\dim_{M}(E)$ of $E$ to be $d_{M}$. 
\end{definition}
\noindent By Proposition 3.4 in \cite{Falconer2004}, in the case
of $M=[0,1]^{d}$, we always have $\dim_{M}(M)=d$, matching the usual
definition of dimension.

We now state the formal assumptions \textbf{\hypertarget{allAssumptions}{A.1-A.6}
}which we will use to prove our main results.
\hluo{These assumptions can be put into three groups, as explained after their statements.}
\label{PageAss1-5}

\begin{enumerate}[wide, labelwidth=!, labelindent=0pt,label={(\roman*)}]
\item Global assumptions on 
\steve{$f$, $M$, and $S_{0}$}

\noindent \textbf{\hypertarget{Assumption1}{A.1.}} (Compact support)
The support of $f$ is $supp(f)=M=[0,1]^{d}$, $d>0$ and we consider
the version of the density $f$ that is zero on $\mathbb{R}^{d}\backslash M$
and continuous on $M$.

\noindent \textbf{\hypertarget{Assumption2}{A.2.}} (Single component)
The zero-density region $S_{0}$ is contained in the interior of $M$
and has one connected component.

\item Local assumptions on $f$ and $S_0$

\noindent \textbf{\hypertarget{Assumption3}{A.3.}} (Order of smoothness)
There exists an $\epsilon_{0}>0$ such that both quantities $\overline{K_{f}}(\epsilon_{0})>0$
and $\underline{K_{f}}(\epsilon_{0})>0$ exist. As
remarked on page \pageref{PageK_eps}, the existence of $\epsilon_{0}>0$ implies that the limiting values $\overline{K_{f}}>0$
and $\underline{K_{f}}>0$ also exist.

\noindent \textbf{\hypertarget{Assumption 4}{A.4.}} Let $\epsilon_{0}>0$
be as described in~A.3. There exist two positive constants $L_{f}$
and $U_{f}$ such that, for all $\epsilon$ with $0<\epsilon<\min(1,\epsilon_{0})$,
$L_{f}\cdot d(\boldsymbol{x},S_{0})^{\overline{K_{f}}(\epsilon)}\leq f(\boldsymbol{x})\leq U_{f}\cdot d(\boldsymbol{x},S_{0})^{\underline{K_{f}}(\epsilon)}$
for all $\boldsymbol{x}\in B_{\epsilon}(S_{0})\cap M$.

\item Assumptions on coverings

\noindent \textbf{\hypertarget{Assumption5}{A.5.}} (Regular covering)
There is an $\eta>0$ such that the elements of the sequence of coverings
$\mathcal{B}_{r(n)}^{d}(M)$ are comprised of balls whose centers
lie in $M$ and whose radii satisfy $r(n)\sim O(n^{-\eta})$. The
sequence of coverings is \textit{regular}, i.e., the cardinalities
$|\mathcal{B}_{r(n)}^{d}(M)|$ of the coverings in the sequence satisfy
the condition $|\mathcal{B}_{r(n)}^{d}(M)|\sim O(n^{d\eta})$.

\noindent \textbf{\hypertarget{Assumption6}{A.6.}} (Restriction
of covering to $S_{0}$) Let $\mathcal{B}_{r(n)}^{d}(M)$ be the covering
considered in A.5. If $d_{0}$ is the Minkowski dimension of $S_{0}$
and $d_{0}<d$, then the number of balls in $\mathcal{B}_{r(n)}^{d}(M)$
intersecting $S_{0}$ is bounded from above by $H_{\varepsilon}(n)\sim O(n^{d_{0}\eta(1+\varepsilon)})$,
for each $\varepsilon>0$. $H_{\varepsilon}(n)$ is a function of sample size $n$ that depends
on the parameter $\varepsilon$.

\end{enumerate}

\hyperlink{allAssumptions}{Assumptions A.1 and A.2} can be regarded
as restrictions on the topological properties of the support of the
density function and we will see later that both assumptions can
be relaxed under suitable conditions. 

\noindent
\hyperlink{allAssumptions}{Assumptions A.3 and A.4} describe the
local behavior of $f$ in the vicinity of $S_{0}$. 
These two assumptions are quite general.  Similar notions appear in both the TDA \cite{guibas2013witnessed} and density estimation literatures \cite{Mammen&Polonik2013}.

\noindent
\hyperlink{allAssumptions}{Assumptions A.5 and A.6} 
are tied
 to the Minkowski dimension of $S_0$ and the covering scheme we
devise for detection of $S_{0}$. 
These two assumptions allow us to describe exactly what we mean by "low dimensionality" of $S_0$.

Next we prove that under \hyperlink{allAssumptions}{Assumptions A.1
and A.2} there exists one covering that satisfies \hyperlink{allAssumptions}{Assumptions
A.5 and A.6} for $S_{0}$ of Minkowski dimension $d_{0}$. 

\begin{lemma}
\label{lem:Exist regular covering} Suppose that A.1 and A.2 hold
and that $S_{0}$ is a lower dimensional zero-density region of Minkowski
dimension $d_{0}<d$. Then there exists a sequence of coverings $\mathcal{B}_{r(n)}^{d}(M)$ that satisfies A.5 and A.6.

\end{lemma}
\begin{proof}
\hluo{
The construction in this lemma is as follows. We first pick a minimal covering of $S_0$ to ensure that the cardinality assumption A.5 is satisfied. Then, we pick a grid-based covering of the support $M$, which ensures that the cardinality assumption A.6 holds. This grid-based covering may introduce more covering balls which intersect with $S_0$. Therefore, we prune our covering, removing the additional covering balls over $S_0$, so that the assumption A.5 still holds. 
This construction requires that $S_0$ be covered by its minimal covering, but covers the rest of the support with grid-based covering balls.
}

\textbf{Step 1.} We prove first that A.6 can be satisfied. \\
Let $r(n) = c n^{-\eta}$, $c > 0$, and let $\mathcal{B}_{r(n)}^{d}(S_0)$ be a minimal $r(n)$ covering of $S_0$ of cardinality $|\mathcal{B}_{r(n)}^{d}(S_0)|$ .
Since we assume that $\dim_{M}(S_{0})=d_{0}$,

\begin{adjustbox}{center}
\parbox{\linewidth}{
\begin{align}
\lim{}_{r(n)\rightarrow0}\frac{\log\mathcal{N}_{r(n)}(S_{0})}{-\log r(n)} =\lim{}_{n\rightarrow\infty}\frac{\log\mathcal{N}_{r(n)}(S_{0})}{\log1/r(n)}=d_{0}\Rightarrow\lim{}_{n\rightarrow\infty}\frac{\log\mathcal{N}_{r(n)}(S_{0})}{\log\left(1/r(n)\right)^{d_{0}}}=1.
\end{align}
}
\end{adjustbox}
Thus, for each $\varepsilon>0$, there exists an $N_{*}(\varepsilon)$ such that, for all $n > N_{*}(\varepsilon)$

\begin{align}
\frac{\log\mathcal{N}_{r(n)}(S_{0})}{\log\left(1/r(n)\right)^{d_{0}}}\leq & 1+\varepsilon, \,\,\,\,\,\, \log\mathcal{N}_{r(n)}(S_{0})\leq  \left(1+\varepsilon\right)\cdot\log\left(1/r(n)\right)^{d_{0}},\nonumber \\ \mathcal{N}_{r(n)}(S_{0})\leq & \left(\frac{1}{r(n)^{d_{0}}}\right)^{1+\varepsilon}
= c^{-d_0 (1 + \varepsilon)} n^{d_0 \eta (1 + \varepsilon)} = H_\varepsilon(n).
\label{eq:key1}
\end{align}

This establishes that each term in the tail of the sequence $\{|\mathcal{B}_{r(n)}^{d}(S_0)|\}$ is bounded by $H_\varepsilon(n)$ as in A.6, 
\steve{ provided that all of the covering balls from the minimal covering are in the final covering constructed in Step 2 below.}

\noindent \textbf{Step 2.} Now we turn to the proof of A.5. 

\textbf{First,} we take $n$ sufficiently large so that there exists
a hyper cube $C(n)$ of dimension $d$ and side length $r(n)$ such
that $B_{3r(n)}(C(n))\subset\text{int }M\backslash B_{5r(n)}(S_{0})$.
For such a fixed $n$ and the $r(n)$-covering $\mathcal{B}_{r(n)}^{d}(S_{0})$
we derived above, we now construct such a covering $\mathcal{B}_{r(n)}^{d,i}(M)$
where the centers of covering balls in $\mathcal{B}_{r(n)}^{d,i}(M)$
are on the grid set
\begin{align*}
G(n)\coloneqq M\cap & \{\boldsymbol{x}=(x_{1},\cdots,x_{d})\in\mathbb{R}^{d}\mid x_{i}\in\mathbb{N}\cdot\frac{r(n)}{d},\forall i=0,\cdots,d\}.  
\end{align*}
For this grid, the maximal distance from an arbitrary point in $M$
to a point in $G(n)$ is at most $\sqrt{d\cdot(\frac{r(n)}{d})^{2}}=\frac{r(n)}{\sqrt{d}}\leq r(n)$
for $d\geq1$. This maximal distance is attained by a pair of points
with coordinates in the form $(x_{1},\cdots,x_{d}),(x_{1}+\frac{r(n)}{d},\cdots,x_{d}+\frac{r(n)}{d})\in G(n)$.
So an arbitrary point in $M$ is contained in some ball of radius
$r(n)$ with its center in $G(n)$. The cardinality of this covering is the same as the cardinality of the set $G(n)$ of centers,
which is of order $O\left(\left(1/\frac{r(n)}{d}\right)^{d}\right)\sim O(n^{\eta d})$.

\textbf{Second,} we delete all covering balls in $\mathcal{B}_{r(n)}^{d,i}(M)$
that intersect $S_{0}$. That is, we obtain another covering $\mathcal{B}_{r(n)}^{d,ii}(M)=\{B\in\mathcal{B}_{r(n)}^{d,i}(M)\mid B\cap S_{0}=\emptyset\}$,
which ensures that there are no additional covering balls intersecting
$S_{0}$ and A.6 is still guaranteed by the covering $\mathcal{B}_{r(n)}^{d}(S_{0})$
we choose in Step 1 because no additional covering balls are touching
$S_{0}$. This operation only deletes those covering balls that are
completely contained in $\mathcal{B}_{3r(n)}^{d}(S_{0})$, so no
covering ball intersecting $C(n)$ will be deleted. As argued
in the first step, $\mathcal{B}_{r(n)}^{d,ii}(M)$ is also a covering
of $C(n)$ of Minkowski dimension $d$ so it must contain covering
balls of cardinality $O(n^{\eta d})$.

\textbf{Third,} we recognize that after deletion and obtaining the covering
$\mathcal{B}_{r(n)}^{d,ii}(M)$, the union $\cup_{B\in\mathcal{B}_{r(n)}^{d,ii}(M)\cup\mathcal{B}_{r(n)}^{d}(S_{0})}B$
may not cover all points of $M$. In the previous deletion operation,
we delete every covering ball such that its center $p$ satisfies
$d(p,S_{0})\leq r(n)$. Therefore any point $q$ in $M\backslash\left(\cup_{B\in\mathcal{B}_{r(n)}^{d,ii}(M)\cup\mathcal{B}_{r(n)}^{d}(S_{0})}B\right)$
will satisfy $d(q,S_{0})\leq2r(n)$. Now we add \emph{finitely many}
translated copies of $\mathcal{B}_{r(n)}^{d}(S_{0})$ to $\mathcal{B}_{r(n)}^{d,ii}(M)\cup\mathcal{B}_{r(n)}^{d}(S_{0})$
to obtain our final covering. Denote by $\mathcal{B}_{r(n)}^{d}(S_{0})+\boldsymbol{v}$
the covering obtained by translating each covering ball in $\mathcal{B}_{r(n)}^{d}(S_{0})$
by vector $\boldsymbol{v}$. We let the translation vector $\boldsymbol{v}$
vary in following set of vectors

\begin{adjustbox}{center}
\parbox{\linewidth}{
\begin{align*}
G_{1}(n) & \coloneqq\{\boldsymbol{v}=(v_{1},\cdots,v_{d})\in\mathbb{R}^{d}\mid v_{i}\in\kappa\cdot\frac{r(n)}{d},\kappa=0,\pm1,\cdots,\pm4d,\forall i=0,\cdots,d\}
\end{align*}
}
\end{adjustbox}
which is a set consisting of $(4d\times2+1)^{d}=(8d+1)^{d}$ vectors.
This construction forms a (irregular) grid $G_{1}(n)$ by translating
the centers of balls in $\mathcal{B}_{r(n)}^{d}(S_{0})$.  This grid
covers all points in $B_{2r(n)}(S_{0})$ since it translates $4r(n)$
in each coordinate direction. The maximal distance from an arbitrary
point in $B_{2r(n)}(S_{0})$ to a point in $G_{1}(n)$ is no more than
$\sqrt{d\cdot(\frac{r(n)}{d})^{2}}=\frac{r(n)}{\sqrt{d}}\leq r(n)$
for $d\geq1$. Therefore 
\begin{align*}
\mathcal{B}_{r(n)}^{d,iii} & \coloneqq\{B\in\mathcal{B}_{r(n)}^{d}(S_{0})+\boldsymbol{v}\mid\boldsymbol{v}\in G_{1}(n)\}
\end{align*}
 will cover the area $B_{2r(n)}(S_{0})$. This will add at most $\left(8d+1\right)^{d}\cdot\left|\mathcal{B}_{r(n)}^{d}(S_{0})\right|\sim O\left(\left|\mathcal{B}_{r(n)}^{d}(S_{0})\right|\right)$
additional covering balls that intersect $S_{0}$. This addition will
not affect the order of magnitude of the cardinality of $\mathcal{B}_{r(n)}^{d,ii}(M)\cup\mathcal{B}_{r(n)}^{d}(S_{0})$
satisfying A.5.  We simply add finitely many $\left(8d+1\right)^{d}$
translated copies of $\mathcal{B}_{r(n)}^{d}(S_{0})$ to cover the
``gaps'' created by deletion in the previous step.

The covering $\mathcal{B}_{r(n)}^{d,\dagger}\coloneqq\mathcal{B}_{r(n)}^{d,ii}(M)\cup\mathcal{B}_{r(n)}^{d}(S_{0})\cup\mathcal{B}_{r(n)}^{d,iii}$
satisfies both A.5 and A.6 by construction as we saw above.
\end{proof}
The constructive proof of the previous lemma makes explicit use of minimal coverings of~$S_0$.  There are many situations in which sequences of coverings satisfying A.5 and A.6 can be produced without resorting to the use of a minimal covering of $S_0$.  The following proposition gives an explicit example of such a situation.
\begin{proposition}
\label{prop:correct-covering-for-line-segment}
For the example in Figure \ref{Minimal Example}, for a fixed $\eta>0$,
there exists a sequence of coverings $\mathcal{B}_{r(n)}^{2}(M)$
satisfying A.5 and A.6.

\begin{figure}[t]
\centering

\includegraphics[scale=0.5]{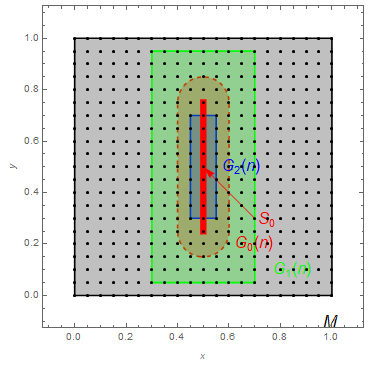}

\caption{Illustration of the construction in Proposition \ref{prop:correct-covering-for-line-segment} with the density in Figure \ref{Minimal Example}. With $r=0.1$,
the grid centers in $G(n)$ are displayed.  The range of $G_{0}(n)$
is displayed in  red; the ranges of $G_{1}(n)$ and $G_{2}(n)$ are
displayed in green and blue, respectively. \label{fig:construction_prop}}
\end{figure}
\end{proposition}

\begin{proof}
We consider the covering where the centers of the balls in $\mathcal{B}_{r(n)}^{2}(M)$
are exactly the grid set
\begin{align}
G(n) & =M\cap\{\boldsymbol{y}=(y_{1},y_{2})\in\mathbb{R}^{2}\mid y_{1},y_{2}\in\mathbb{N}\cdot\frac{r(n)}{2}\}.
\end{align}
(See Figure~\ref{fig:construction_prop} for an illustration.)
The cardinality of this covering is the same as the cardinality
of the grid set $G(n)$. As $n$ increases, we have a sequence of
coverings whose cardinalities are of order $O\left(\frac{1}{r(n)/2}\right)^{2}\sim O(n^{2\eta})$
as required in A.5. 

Next, we prove that the covering balls of $\mathcal{B}_{r(n)}^{2}(M)$
intersecting $S_{0}$ also have cardinalities that satisfy A.6. Let
$n$ be large enough that $r(n)<1/10$ and the $r(n)$-neighborhood
of $S_{0}$, $B_{r(n)}(S_{0})\subset\text{int }M$. Consider the set
of grids $G_{0}(n)\coloneqq\{\boldsymbol{y}=(y_{1},y_{2})\in G(n)\mid0\leq d(\boldsymbol{y},S_{0})<r(n)\}$.
A covering ball in $\mathcal{B}_{r(n)}^{2}(M)$ intersects $S_{0}$
if and only if its center is in $G_{0}(n)$. We construct the rectangular
grids defined below to bound the cardinality of $G_{0}(n)$.
\begin{adjustbox}{center}
\parbox{\linewidth}{
\begin{align*}
G_{1}(n)\coloneqq & \{\boldsymbol{y}=(y_{1},y_{2})\in G(n)\mid\frac{1}{2}-2r(n)\leq y_{1}\leq\frac{1}{2}+2r(n),\frac{1}{4}-2r(n)\leq y_{2}\leq\frac{3}{4}+2r(n)\},\\
G_{2}(n)\coloneqq & \{\boldsymbol{y}=(y_{1},y_{2})\in G(n)\mid\frac{1}{2}-\frac{r(n)}{2}\leq y_{1}\leq\frac{1}{2}+\frac{r(n)}{2},\frac{1}{4}+\frac{r(n)}{2}\leq y_{2}\leq\frac{3}{4}-\frac{r(n)}{2}\}.
\end{align*}
}
\end{adjustbox}
It is clear from their definitions that $G_{2}(n)\subset G_{0}(n)\subset G_{1}(n)$.

In the grid set $G_{2}(n)$, there are at least $\frac{r(n)}{r(n)/2}=2$
grid points with the same value of $y_{2}$; and at least $\left\lfloor \frac{1/2-r(n)}{r(n)/2}\right\rfloor =\left\lfloor \frac{1}{r(n)}-2\right\rfloor $
points with the same value of $y_{1}$. Thus there are at least $2\cdot\left\lfloor \frac{1}{r(n)}-2\right\rfloor $
grid points in $G_{2}(n)$. Similarly, in the grid set $G_{1}(n)$
there are at most $\frac{4r(n)}{r(n)/2}+1=9$ grid points with the
same value of $y_{2}$; and at most $\left\lceil \frac{1/2+4r(n)}{r(n)/2}+1\right\rceil =\left\lceil \frac{1}{r(n)}+9\right\rceil $
points with the same value of $y_{1}$. Thus there are at most $9\cdot\left\lceil \frac{1}{r(n)}+9\right\rceil $
grid points in $G_{2}(n)$. 

Since $G_{2}(n)\subset G_{0}(n)\subset G_{1}(n)$, the number of balls
intersecting $S_{0}$ has cardinality $\left|G_{0}(n)\right|$ bounded
by the inequality 
\begin{equation}
2\cdot\left\lfloor \frac{1}{r(n)}-2\right\rfloor \leq\left|G_{2}(n)\right|\leq\left|G_{0}(n)\right|\leq\left|G_{1}(n)\right|\leq9\cdot\left\lceil \frac{1}{r(n)}+9\right\rceil .\label{eq:key_bounding_ineq}
\end{equation}
Both sides of (\ref{eq:key_bounding_ineq}) have order of magnitude
$O\left(\frac{1}{r(n)}\right)\sim O(n^{\eta})$ and this verifies
A.6.
\end{proof}

Although we have shown by example that our \hyperlink{allAssumptions}{Assumption
A.1-A.6} can be verified for many densities,
we point out that in general it might not be easy to verify
these assumptions when $S_{0}$ or $M$ possesses a complicated topological
structure.

\subsection{Main result }
\hluo{Our} main result connects coverings to data \hluo{points}.  An ``empty'' ball centered at $\boldsymbol{x}$ is one for which $\{X_{1},\cdots,X_{n}\}\cap B_{r(n)}(\boldsymbol{x})=\emptyset$. 
\begin{theorem}
\label{Main Theorem} Suppose we have a set of i.i.d.$\,$data $X_{1},\cdots,X_{n}$
drawn from a continuous density $f(\boldsymbol{x})$ w.r.t.$\,$the
$d$-dimensional Lebesgue measure $\nu^{d}$ defined on a compact
subset $M=[0,1]^{d}$ of $\mathbb{R}^{d}$ and assume that \hyperlink{allAssumptions}{Assumptions
A.1-A.6} hold. Assume also that the radius $r$ and the separation distance $\epsilon$
satisfy the following growth rates:
\begin{align*}
r(n)\sim O(n^{-\eta}) & ,\,0<\eta<\frac{1}{d},\\
\epsilon(n)\sim O(n^{-\psi}) & ,\:0<\psi\leq\eta,\\
2r & (n)\leq\epsilon(n)<1.
\end{align*}
Finally, assume the validity of the following bounding condition for the density $f$ outside the
$\epsilon$-neighborhood of $S_{0}$: 
\begin{align*}
m(f,n)\coloneqq\min_{\boldsymbol{w}\in M\backslash B_{\epsilon}(S_{0})}f(\boldsymbol{w}) & \in(0,\infty)\sim O(n^{-\xi}),\,0<\xi<\frac{1-2\eta d}{2}.
\end{align*}
Then:\\
(A) If $\eta$ and  $\psi$ satisfy 
\begin{alignat*}{1}
\begin{array}{c}
1-2\eta d-2\overline{K_{f}}\psi>0,\end{array}
\end{alignat*}
we have $\lim_{n\rightarrow\infty}P(\text{\hluo{there are no empty} \ensuremath{\epsilon}(n)-outside balls})=1$.

\noindent 
(B) If $\eta$ satisfies 
\begin{align*}
\begin{array}{c}
1+d_{0}\eta-\underline{K_{f}}\eta-d\eta<0,\end{array}
\end{align*}
we have $\lim_{n\rightarrow\infty}P(\text{all \ensuremath{\epsilon}(n)-inside balls are empty})=1$. 
\end{theorem}

\hluo{\noindent The proof is given in Appendix \ref{Proof of the Main Theorem}.}

Similarly, we have a result dealing with the situation where $S_0$ has more than
one (but a finite number) connected component. 
\begin{corollary}
\label{cor:Suppose-we-have}
Under the same assumptions as in Theorem
\ref{Main Theorem}, 
suppose that the zero-density
region 
can be decomposed into $K$ disjoint
connected components,
$S_{0}=\sqcup_{k=1}^{K}S_{0k}$. 
For $n$ large enough (so that $\epsilon$ is small enough),
an $\epsilon$-neighborhood of $S_0$ will also be comprised of
$K$ disjoint
connected components, the $\epsilon$-neighborhoods of each component $S_{0k}$. 

Assume that the following
bounding conditions for the density
$f$ outside the $\epsilon$-neighborhood of each component $S_{0k}$ hold:
\begin{align*}
m_{k}(f,n)\coloneqq\min_{\boldsymbol{w}\in M\backslash B_{\epsilon(n)}(S_{0k})}f(\boldsymbol{w}) & \in(0,\infty)\sim O(n^{-\xi}),\,0<\xi<\frac{1-2\eta d}{2}.
\end{align*}
Assume also that the orders of smoothness of $f$ w.r.t.~$S_{0k}$ are $\overline{K_{f,k}}$ and $\underline{K_{f,k}}>0$
and that, for $S_{0k}$ of dimension $d_{0k}$, the cardinality of the set of covering
balls intersecting $S_{0k}$ is of order $O(n^{\eta d_{0k}})$.

\noindent
Then:

\noindent 
(A) If $\eta$ and $\psi$ satisfy 
\begin{alignat*}{1}
\min_{k=1,\cdots,K}\begin{array}{c}
1-2\eta d-2\overline{K_{f,k}}\psi>0,\end{array}
\end{alignat*}
we have $\lim_{n\rightarrow\infty}P(\text{\hluo{there are no empty} \ensuremath{\epsilon}(n)-outside balls})=1$
for all $k$.

\noindent
(B) If $\eta$ satisfies
\begin{align*}
\max_{k=1,\cdots,K}\begin{array}{c}
1+d_{0k}\eta-\underline{K_{f,k}}\eta-d\eta<0,\end{array}
\end{align*}
we have $\lim_{n\rightarrow\infty}P(\text{all \ensuremath{\epsilon}(n)-inside balls are empty})=1$
for all $k$. 
\end{corollary}

\section{\label{sec:Non-compact-support}Non-compact support}

The main result of the previous section relies on the assumption of a
compact support for the density $f$. However, in practice, many distributions
have unbounded supports. 
In this section we consider the case of a density function $f$ with a non-compact support contained in $\mathbb{R}^d$ and, assuming that the 
tails of the density decay at certain rates, we derive 
results similar to those of the previous section.

Our strategy restricts our consideration to a region within the support
that contains most of the probability mass of the density $f$. The
original motivation for these ideas comes from the concentration inequalities
that describe the fact that observations usually ``concentrate''
around the ``center'' of a probability density with high probability
(for example, \cite{Talagrand2014}).

Suppose that the support of the probability density function $f$
is a non-compact set $M^{'}\subset\mathbb{R}^{d}$.  We consider a
truncation of the non-compact support $M^{'}$ to a compact subset
that contains most of the probability mass
of $f$ and whose interior contains (all components of) $S_{0}$.
This allows us to examine the topological features of $M^{'}$ over a compact
truncated region that contains the bulk of the data and to establish our results
using
arguments similar to those used to prove Theorem~\ref{Main Theorem}. 
For ease of exposition, we focus parts of our discussion on the situation where $M^{'}=\mathbb{R}^{d}$, although most of our results 
can be extended to more general situations. We will state the restrictions on the 
tail behavior of the
density functions in Proposition~\ref{prop:Polynomial-Exp tail}.
\begin{definition}
\label{1-delta support}($(1-\delta)$-support containing $S_{0}$)
A \emph{$(1-\delta)$-support }containing \emph{$S_{0}$} for some
$0<\delta<1$, denoted by $S_{f}^{1-\delta}\subset\mathbb{R}^{d}$,
is a subset of the form $S_{f}^{1-\delta}\coloneqq[-B,B]^{d}\subset\mathbb{R}^{d}$
for some $B\in(0,\infty)$, such that $S_{0}\subset\text{int }S_{f}^{1-\delta}=(-B,B)^{d}$
and $P\left\{ X\in S_{f}^{1-\delta}\right\} =1-\delta$. 
\end{definition}
The hypercube $[-B,B]^{d}$ containing $1-\delta$ of the probability
expands as $\delta$ shrinks. If $S_{0}$ is compact then $S_{f}^{1-\delta}$
exists for a small enough $\delta_{0}\in(0,1)$; if $S_{0}$ is non-compact,
then it is unbounded and $S_{f}^{1-\delta}$ does not exist for any
$\delta\in(0,1)$. For our subsequent results, we assume that the $(1-\delta)$-support
is completely contained in the interior of the support $M^{'}$.

Passing from $M^{'}$ to $S_{f}^{1-\delta}$ allows us to work with
a compact cubical subset and removes the technical difficulties associated with densities whose tails decay to zero. We replicate the results
of Section \ref{sec:Statement-of-the-main-result} with modifications
to \hyperlink{allAssumptions}{Assumptions A.1-A.6}. The modifications
stem from replacing the compact support $[0,1]^{d}$ with a compact
set $S_{f}^{1-\delta}$ that covers $S_{0}$ and contains $1-\delta$
mass of the probability measure. 

We first consider the case where $\delta>0$ and $S_{f}^{1-\delta}$
are both fixed and state \hyperlink{allAssumptions}{Assumptions
A.1'-A.6'} for the non-compact support case.

\begin{enumerate}[wide, labelwidth=!, labelindent=0pt,label={(\roman*)}]
\item \steve{Global assumptions on $f$, $M$, and $S_{0}$}

\noindent \textbf{\hypertarget{Assumption1}{A.1'.}} (Non-compact
support) The support of $f$ is $supp(f)=M^{'}\subset\mathbb{R}^{d}$,
$d>0$ and we consider the version of the density $f$ that is zero
on $\mathbb{R}^{d}\backslash M^{'}$ and continuous on $M^{'}$.

\noindent \textbf{\hypertarget{Assumption2}{A.2'.}} (Single component)
The zero-density region $S_{0}$ is contained in the interior of $M^{'}$
and has one connected component.

\item \steve{Local assumptions on $f$ and $S_{0}$}

\noindent \textbf{\hypertarget{Assumption3}{A.3'.}} (Order of smoothness)
There exists an $\epsilon_{0}>0$ such that both quantities $\overline{K_{f}}(\epsilon_{0})>0$
and $\underline{K_{f}}(\epsilon_{0})>0$ exist. 

\noindent \textbf{\hypertarget{Assumption 4}{A.4'.}} Let $\epsilon_{0}>0$
be as described in~A.3'. There exist two positive constants $L_{f}$
and $U_{f}$ such that, for all $\epsilon$ with $0<\epsilon<\min(1,\epsilon_{0})$,
$L_{f}\cdot d(\boldsymbol{x},S_{0})^{\overline{K_{f}}(\epsilon)}\leq f(\boldsymbol{x})\leq U_{f}\cdot d(\boldsymbol{x},S_{0})^{\underline{K_{f}}(\epsilon)}$
for all $\boldsymbol{x}\in B_{\epsilon}(S_{0})\cap S_{f}^{1-\delta}$.

\item \steve{Assumptions on coverings}

\noindent \textbf{\hypertarget{Assumption5}{A.5'.}} (Regular covering)
There is an $\eta>0$ such that the elements of the sequence of coverings
$\mathcal{B}_{r(n)}^{d}(S_{f}^{1-\delta})$ are comprised of balls
whose centers lie in $\text{int }S_{f}^{1-\delta}$\textbf{\textcolor{red}{{}
}}and whose radii satisfy $r(n)\sim O(n^{-\eta})$. The sequence of
coverings is \textit{regular}, i.e., the cardinalities $|\mathcal{B}_{r(n)}^{d}(S_{f}^{1-\delta})|$
of the coverings in the sequence satisfy the condition $|\mathcal{B}_{r(n)}^{d}(S_{f}^{1-\delta})|\sim O(n^{d\eta})$.

\noindent \textbf{\hypertarget{Assumption6}{A.6'.}} (Restriction
of covering to $S_{0}$) Let $\mathcal{B}_{r(n)}^{d}(S_{f}^{1-\delta})$
be the covering considered in A.5'. If $d_{0}$ is the Minkowski dimension
of $S_{0}$ and $d_{0}<d$, then the number of balls in $\mathcal{B}_{r(n)}^{d}(S_{f}^{1-\delta})$
intersecting $S_{0}$ is bounded from above by $H_{\varepsilon}(n)\sim O(n^{d_{0}\eta(1+\varepsilon)})$,
for each $\varepsilon>0$.

\end{enumerate}

As for the compact case on page~\pageref{PageAss1-5}, these assumptions can be divided into three groups. A.1' and A.2' are restrictions on the topology on the support; A.3' and A.4' are descriptions of the local behavior of the density function; and A.5' and A.6' stipulate the existence of a sequence of coverings with good properties.
\begin{theorem}
\label{Noncompact Theorem} Consider a sequence of i.i.d.$\,$data
$X_{1},X_{2}\cdots$ drawn from a distribution having a continuous density $f(\boldsymbol{x})$
w.r.t.$\,$the $d$-dimensional Lebesgue measure $\nu^{d}$ 
on $\mathbb{R}^{d}$. Assume that \hyperlink{Assumption1new}{A.1'-A.6'}
hold and that \textcolor{black}{$S_{f}^{1-\delta}\subset\text{int }M^{'}$.}
Assume also that the radius $r$ and the separation distance $\epsilon$
satisfy the following growth rates:
\begin{align*}
r(n)\sim O(n^{-\eta}) & ,\,0<\eta<\frac{1}{d},\\
\epsilon(n)\sim O(n^{-\psi}) & ,\,0<\psi\leq\eta,\\
2r & (n)\leq\epsilon(n)<1.
\end{align*}
Finally, assume that the density $f$ outside the $\epsilon$-neighborhood of
$S_{0}$ satisfies: 
\begin{align*}
m(f,n)\coloneqq\min_{\boldsymbol{w}\in S_{f}^{1-\delta}\backslash B_{\epsilon}(S_{0})}f(\boldsymbol{w}) & \in(0,\infty)\sim O(n^{-\xi}),\,0<\xi<\frac{1-2\eta d}{2}.
\end{align*}
Then:

\noindent 
(A) If $\eta$ and $\psi$ satisfy 
\begin{alignat*}{1}
\begin{array}{c}
1-2\eta d-2\overline{K_{f}}\psi>0,\end{array}
\end{alignat*}
we have $\lim_{n\rightarrow\infty}P(\text{\hluo{there are no empty }\ensuremath{\epsilon}(n)-outside balls})=1$.

\noindent 
(B) If $\eta$ and $\psi$ satisfy 
\begin{align*}
\begin{array}{c}
1+d_{0}\eta-\underline{K_{f}}\eta-d\eta<0,\end{array}
\end{align*}
we have $\lim_{n\rightarrow\infty}P(\text{all \ensuremath{\epsilon}(n)-inside balls are empty})=1$. 
\end{theorem}

\begin{proof}
The idea of the proof is that the number of samples falling inside
$S_{f}^{1-\delta}$ is close to the total number of samples, $n$,
by the definition of $(1-\delta)$-support and the law of large numbers.  By assumption, we know that
$S_{f}^{1-\delta}$ is the hyper-cube $[-B,B]^{d}$ contained in $\text{int }M^{'}$
and that \hyperlink{Assumption1}{A.1'-A.6'} hold. The number $n_{\delta}$ of observations that fall in $S_{f}^{1-\delta}$
has a $Binom(n,1-\delta)$ distribution.
Using the two-sided Hoeffding inequality, (\ref{eq:Hoeffding}) of the appendix, with
$\gamma=\frac{\delta}{2},p=1-\delta$, we have:{\footnotesize{} }
\begin{align*}
P(n_{\delta}\geq n(1-\delta-\frac{\delta}{2})) & =P(n_{\delta}\geq n(1-\frac{3\delta}{2}))\\
 & \geq P((1-\frac{3\delta}{2})n\leq n_{\delta}\leq(1-\frac{\delta}{2})n)\\
 & \geq1-2\exp\left(-\frac{\delta}{2}^{2}n\right).
\end{align*}
Conditioning on the event that that $n_{\delta}\geq n(1-\frac{3\delta}{2})$
merely adds a nonzero multiplicative factor to the probabilities of
the events appearing in cases (A) and (B) of the statement of the theorem. With this adjustment
the proof of this theorem parallels that of Theorem \ref{Main Theorem}.
The conditional probabilities of the events in (A) and (B) tend to 1 as $n\rightarrow\infty$.
We also know that, as $n\rightarrow\infty$, the probability of the
conditioning event tends to 1, completing the proof. 
\end{proof}
The result above is stated for a fixed $(1-\delta)$-support $S_{f}^{1-\delta}$. As the sample size $n$ tends to infinity, we do not want to constrain
ourselves to a fixed $S_{f}^{1-\delta}$, since this will prevent us
from exploring the entire support of $f$. Consider a sequence of
regions $S_{f}^{1-\delta(n)}$ with a decreasing sequence of $\delta=\delta(n)$,
expanding the region under consideration as the sample accumulates.
To retain the conclusions of the theorem, we must take care to expand
the region slowly enough to control the decay of the density near
the edges of $S_{f}^{1-\delta(n)}$.  

In what follows, we denote by $\mathcal{B}_{r(n)}^{d}$ a covering
of $S_{f}^{1-\delta(n)}$ and assume that A.5' and A.6' hold for this sequence of coverings.
(Note that, here, $\mathcal{B}_{r(n)}^{d}$ denotes a covering of the $(1-\delta)$ support, not of $\mathbb{R}^d$, and that we dropped the explicit dependence on $S_{f}^{1-\delta(n)}$
that appears in A.5' and A.6' to simplify notation.)

Theorem \ref{Noncompact Theorem} considers a fixed $\delta$.  The upcoming theorems \ref{Decreasing sequence of deltas} and \ref{Decreasing sequence of deltas-1} allow $\delta$ to shrink toward $0$. 
\begin{theorem}
\label{Decreasing sequence of deltas}Under the assumptions of Theorem
\ref{Noncompact Theorem} 
(except for $\delta(n)$)
, consider a positive decreasing sequence
$\delta(n)\,\text{satisfying }\lim_{n\rightarrow\infty}\delta(n)=0$
and the associated sequence of $(1-\delta)$-supports $S_{f}^{1-\delta(n)}$.
For each $(1-\delta)$-support $S_{f}^{1-\delta(n)}$ we consider
its corresponding covering $\mathcal{B}_{r(n)}^{d}$ consisting of
open balls of radius $r(n)$.

Assume that 
\begin{align*}
M(f)\coloneqq\sup_{\boldsymbol{w}\in\mathbb{R}^{d}}f(\boldsymbol{w}) & \in(0,\infty),\\
m(f,n)\coloneqq\min_{\boldsymbol{w}\in S_{f}^{1-\delta(n)}\backslash B_{\epsilon(n)}(S_{0})}f(\boldsymbol{w}) & \in(0,\infty)\sim O(n^{-\xi}),\,0<\xi<\frac{1-2\eta d}{2},
\end{align*}
and that the cardinality $|\mathcal{B}_{r(n)}^{d}|\sim o(n^{\Omega}),\,\Omega\geq\eta>0$. 

\noindent
Then, if $\eta$ and $\psi$ satisfy 
\begin{alignat*}{1}
\begin{array}{c}
1-2\eta d-2\overline{K_{f}}\cdot\psi>0,\end{array}
\end{alignat*}
we have $\lim_{n\rightarrow\infty}P(\text{\hluo{there are no empty }\ensuremath{\epsilon}(n)-outside balls})=1$.
\end{theorem}

\begin{proof}
We denote the number of total samples by $n$ and the number of samples
falling inside $S_{f}^{1-\delta(n)}$ by $\mathcal{N}_{\delta(n)}$.
To start, let us consider a sufficiently large $n\geq\mathsf{N}_{1}$
such that the decreasing sequence satisfies $\delta(n)\leq0.1$ so that $1-\delta(n)\ge0.9$.
Then, let us decompose the event $C_{n}=\{\text{\hluo{there are no empty }\ensuremath{\epsilon(n)}-outside balls}\}$
using the event $D_{n}=\{0.95n\leq\mathcal{N}_{\delta(n)}\leq n\}$
and its complement. The the law of total probability yields: 
\begin{align}
P(C_{n}) & =P(C_{n}\mid D_{n})\cdot P(D_{n})+P(C_{n}\mid D_{n}^{c})\cdot P(D_{n}^{c})\nonumber \\
 & \geq P(C_{n}\mid D_{n})\cdot P(D_{n}).\label{eq:decompose bound}
\end{align}

The event $D_{n}$ involves a binomial random variable $\mathcal{N}_{\delta(n)}$,
therefore $P(D_{n})$ can be bounded by the Hoeffding inequality, (\ref{eq:Hoeffding}) of the appendix.
Since $1-\delta(n)\ge0.9$, 
\begin{align}
P(D_{n})=P(0.95n & \leq\mathcal{N}_{\delta(n)}\leq n)\geq1-2\exp\left(-\frac{0.1}{2}^{2}n\right).\label{eq:est_n_delta}
\end{align}

The other term, $P(C_{n}\mid D_{n})$, can also be bounded by the
following argument. (Note that $r(n)$ and $\epsilon(n)$ only depend
on the total sample size $n$, not on $\mathcal{N}_{\delta(n)}$.)
\begin{align}
P(C_{n}\mid D_{n}) & =1-P(C_{n}^{c}\mid D_{n})\nonumber \\
 & \geq1-\sum_{B\in\{\epsilon(n)-\text{outside balls}\}}P(B\text{ empty}\mid D_{n})\nonumber \\
 & \geq1-|\mathcal{B}_{r(n)}^{d}|\cdot\left(1-\inf_{B\in\{\epsilon(n)-\text{outside balls}\}}P(B\text{ not empty}\mid D_{n})\right).\label{eq:key to non-compact-2-1}
\end{align}

By the bound in (\ref{epsilon-outside ball bound-1}) in the proof
of Theorem \ref{Main Theorem} in the appendix and the assumption that $|\mathcal{B}_{r(n)}^{d}|\sim o(n^{\Omega})$
we have
\[
|\mathcal{B}_{r(n)}^{d}|\cdot\left(1-\inf_{B\in\{\epsilon(n)-\text{outside balls}\}}P(B\text{ not empty}\mid D_{n})\right)\sim
\]
\[
o(n^{\Omega})\cdot O\left(\exp\left(-n^{\min(1-2\eta d-2\overline{K_{f}}\cdot\psi,1-2\eta d-2\xi)}\right)\right),
\]
which is monotonically decreasing since $\min(1-2\eta d-2\overline{K_{f}}\cdot\psi$ \steve{)} $> 0$ and $\xi<\frac{1-2\eta d}{2}$ by the assumptions of the theorem.
Using the definition of big/small
O notation, we can find a constant $L\in(0,\infty)$ such that
\begin{align}
P(C_{n}\mid D_{n}) & \geq1-L\cdot n^{\Omega}\cdot\exp\left(-n^{1-2\eta d-2\overline{K_{f}}\cdot\psi}\right)\text{ for sufficiently large }n.\label{eq:key to non-compact-2-2}
\end{align}
The bounds (\ref{eq:est_n_delta}) and (\ref{eq:key to non-compact-2-2})
can be substituted back into (\ref{eq:decompose bound}) to obtain 
\begin{align}
P(C_{n}) & \geq P(C_{n}\mid D_{n})\cdot P(D_{n})\nonumber \\
 & \geq\left(1-L\cdot n^{\Omega}\cdot\exp\left(-n^{1-2\eta d-2\overline{K_{f}}\cdot\psi}\right)\right)\left(1-2\exp\left(-\frac{0.1}{2}^{2}n\right)\right).
\end{align}
As long as 
$1-2\eta d-2\overline{K_{f}}\cdot\psi>0$,
this lower bound tends to 1 as $n$ goes to infinity.
\end{proof}
It is important to notice that it is the decay rates of $m(f,n)$ and
$|\mathcal{B}_{r(n)}^{d}|\sim o(n^{\Omega})$ that determine the sufficient
conditions of Theorem \ref{Decreasing sequence of deltas}, not the
rate at which the volume of $S_{f}^{1-\delta(n)}$ increases.

The following theorem deals with the $\epsilon(n)$-inside balls.
\begin{theorem}
\label{Decreasing sequence of deltas-1}Under the assumptions of Theorem
\ref{Noncompact Theorem} 
(except for $\delta(n)$)
, consider a positive decreasing sequence
$\delta(n)\,\text{satisfying }\lim_{n\rightarrow\infty}\delta(n)=0$
and the associated sequence of $(1-\delta)$-supports $S_{f}^{1-\delta(n)}$.
For each $(1-\delta)$-support $S_{f}^{1-\delta(n)}$ we consider
its corresponding covering $\mathcal{B}_{r(n)}^{d}$ consisting of
open balls of radius $r(n)$.

Assume that 
\begin{align*}
m(f,n)\coloneqq\min_{\boldsymbol{w}\in S_{f}^{1-\delta(n)}\backslash B_{\epsilon(n)}(S_{0})}f(\boldsymbol{w}) & \in(0,\infty)\sim O(n^{-\xi}),\,0<\xi<\frac{1-2\eta d}{2}.
\end{align*}

\noindent 
Then, if $\eta$ and $\psi$ satisfy 
\begin{align*}
\begin{array}{c}
1+d_{0}\eta-\underline{K_{f}}\eta-d\eta<0,\end{array}
\end{align*}
we have $\lim_{n\rightarrow\infty}P(\text{all \ensuremath{\epsilon}(n)-inside balls are empty})=1$. 
\end{theorem}

\begin{figure}[t]
\centering\includegraphics[width=7cm,height=4.5cm]{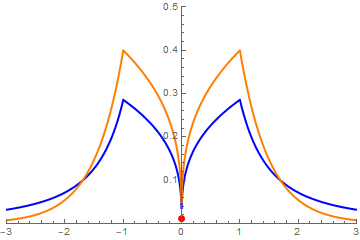}\caption{Density functions with polynomial tails (blue) and exponential tails (orange).\label{fig:tails}}
\end{figure}

\begin{proof}
This result can be established by following exactly the same proof as that for case (B) in Theorem~\ref{Noncompact Theorem}.
Since $S_{0}\subset\text{int }M'$ we can choose $\mathsf{N}_{1}$
sufficiently large so that $B_{\epsilon(n)}(S_{0})\subset\text{int }S_{f}^{1-\delta(n)}$,
for $n\geq\mathsf{N}_{1}$. The probability of the event that an individual
$\epsilon(n)$-inside ball is empty does not depend on the varying
sequence of regions $S_{f}^{1-\delta(n)}$. For 
an appropriate $\mathsf{N}_{2}$ and
every $n\geq\mathsf{N}_{2}$,
our assumption A.6' restricts the number of $\epsilon$-inside balls and
this cardinality does not depend on the varying $(1-\delta(n))$-supports.
The rest of the proof parallels that of Theorem \ref{Noncompact Theorem}
(B), considering $n\geq\max(\mathsf{N}_{1},\mathsf{N}_{2})+1$.
\end{proof}
The existence of a sequence of $\delta(n)$ values satisfying the above conditions
can be verified for densities exhibiting specific tail behaviors on their supports.  
Two examples of such densities are 
given in the two examples 
following Definitions~\ref{def:poly} and~\ref{def:expo}
and are depicted in Figure~\ref{fig:tails}.

\begin{definition}\label{def:poly}
(Polynomial tail) We say that a continuous density $f$ supported
on $\mathbb{R}^{d}$ has a \emph{polynomial tail} if it has the form
\begin{equation}
f(\boldsymbol{x})=\begin{cases}
f_{1}(d(\boldsymbol{x},S_{0}))=C_{1}\cdot d(\boldsymbol{x},S_{0})^{\gamma}, & d(\boldsymbol{x},S_{0})<\epsilon_{0},\\
f_{2}(d(\boldsymbol{x},S_{0}))=C_{2}\cdot d(\boldsymbol{x},S_{0})^{\chi}, & d(\boldsymbol{x},S_{0})\geq\epsilon_{0},
\end{cases}\label{eq:poly tail eq}
\end{equation}
 for some $C_{1},C_{2}\in(0,\infty)$, $\gamma>0$, and $\chi<-d$.
Note that the continuity assumption on the density $f$ requires $C_{1}\cdot\epsilon_{0}^{\gamma}=C_{2}\cdot\epsilon_{0}^{\chi}$.
\end{definition}

\begin{example}
\label{ex_poly}
Consider a continuous density $f$ of the form (\ref{eq:poly tail eq})
with $d=1,\,S_{0}=\{0\},\,\epsilon_{0}=1,\,\gamma=\frac{1}{3},\,\text{and }\chi=-2$,
leading to the density $$f(x)=\begin{cases}
\frac{2}{7}\cdot|x|^{1/3}, & |x|<1,\\
\frac{2}{7}\cdot|x|^{-2}, & |x|\geq1,
\end{cases}.$$ This density yields $\frac{2}{7}\int_{-1}^{1}|x|^{1/3}dx=\frac{3}{7}$.
We can take $n$ sufficiently large so that $\delta(n)<\frac{4}{7}$.
From the definition of $(1-\delta(n))$-support we have $B(n)>\epsilon_{0}=1$
for $\delta(n)<\frac{4}{7}$. The defining equation of $B(n)$ is $\int_{-\infty}^{-B(n)}f(x)dx+\int_{B(n)}^{\infty}f(x)dx=2\int_{B(n)}^{\infty}f(x)dx=\delta(n).$
Therefore $B(n)=\frac{1}{\frac{7}{4}\delta(n)}\sim O(\delta(n)^{-1})$. 

The maximum of this density function is attained at $x=\pm1$ and
equals $M(f)=\frac{2}{7}$. Let $\xi$ and $\psi$ be as defined in
Theorem \ref{Noncompact Theorem}. Then, if we choose sequences $\delta(n)\sim O(n^{\xi/\chi})\sim O(n^{-\xi/2})$
and $\text{\ensuremath{\epsilon(n)}}\sim O(n^{-\psi})\sim O(n^{-\xi/\gamma})\sim O(n^{-3\xi})$,
we will ensure that the minimal value $m(f,n)\coloneqq\min_{\boldsymbol{w}\in S_{f}^{1-\delta(n)}\backslash B_{\epsilon(n)}(S_{0})}$ $f(\boldsymbol{w})=$ $\min\left(C_{1}\cdot\epsilon(n)^{\gamma},C_{2}\cdot B(n)^{\chi}\right)=\min\left(\frac{2}{7}\epsilon(n)^{1/3},\frac{2}{7}B(n)^{-2}\right)$.
We have $\frac{2}{7}\epsilon(n)^{1/3}$ $\sim O(n^{-\xi})$ and $\frac{2}{7}B(n)^{-2}\sim O(n^{-\xi})$,
so that $m(f,n)\sim O(n^{-\xi})$. Furthermore, $\delta(n)\sim O(n^{-\xi/2})\rightarrow0$.
The assumptions on the density function in Theorem \ref{Noncompact Theorem}
are therefore satisfied.
\end{example}

\begin{definition}\label{def:expo}
(Exponential tail) We say that a continuous density $f$ supported
on $\mathbb{R}^{d}$ has an \emph{exponential tail} if it has the
form 
\begin{equation}
f(\boldsymbol{x})=\begin{cases}
f_{1}(d(\boldsymbol{x},S_{0}))=C_{1}\cdot d(\boldsymbol{x},S_{0})^{\gamma}, & d(\boldsymbol{x},S_{0})<\epsilon_{0},\\
f_{2}(d(\boldsymbol{x},S_{0}))=C_{2}\cdot\exp(\beta d(\boldsymbol{x},S_{0})), & d(\boldsymbol{x},S_{0})\geq\epsilon_{0},
\end{cases}\label{eq:exp tail eq}
\end{equation}
 for some $C_{1},C_{2}\in(0,\infty)$, $0<\gamma$ and $\beta<0$.
Note that the continuity assumption on the density $f$ requires $C_{1}\cdot\epsilon_{0}^{\gamma}=C_{2}\cdot\exp(-\beta\epsilon_{0})$. 
\end{definition}

\begin{example}
\label{ex_expo}
Consider a continuous density $f$ of the form (\ref{eq:exp tail eq})
with $d=1,\,S_{0}=\{0\},\,\epsilon_{0}=1,\,\gamma=\frac{1}{3},\,\text{and }\beta=-2$,
leading to the density $$f(x)=\begin{cases}
\frac{2}{5}\cdot|x|^{1/3}, & |x|<1,\\
\frac{2}{5}\cdot e^{2}\cdot e^{-2|x|}, & |x|\geq1,
\end{cases}.$$ This density yields $\int_{-1}^{1}\frac{2}{5}\cdot|x|^{1/3}dx=\frac{3}{5}$.
We can take $n$ sufficiently large so that $\delta(n)<\frac{2}{5}$.
From the definition of $(1-\delta(n))$-support we have $B(n)>\epsilon_{0}=1$
for $\delta(n)<\frac{2}{5}$. The defining equation of $B(n)$ is
$\int_{-\infty}^{-B(n)}f(x)dx+\int_{B(n)}^{\infty}f(x)dx=2\int_{B(n)}^{\infty}f(x)dx=\delta(n).$
Therefore $B(n)=1-\frac{1}{2}\log\left(1+\frac{5}{2}\delta(n)\right)\sim O(-\log\delta(n))$. 

The maximum of this density function is attained at $x=\pm1$ and
equals $M(f)=\frac{2}{5}$. Let $\xi$ and $\psi$ be as defined in
Theorem \ref{Noncompact Theorem}. Then if we choose sequences $\delta(n)\sim O(n^{\xi/\beta})\sim O(n^{-\xi/2})$
and $\text{\ensuremath{\epsilon(n)}}\sim O(n^{-\psi})\sim O(n^{-\xi/\gamma})\sim O(n^{-3\xi})$,
we will ensure that the minimal value $m(f,n)\coloneqq\min_{\boldsymbol{w}\in S_{f}^{1-\delta(n)}\backslash B_{\epsilon(n)}(S_{0})}f(\boldsymbol{w})=\min\left(C_{1}\cdot\epsilon(n)^{\gamma},C_{2}\cdot\exp(\beta B(n))\right)$
\newline
$=\min\left(\frac{2}{5}\epsilon(n)^{1/3},\frac{2}{5}e^{2}\cdot\exp\left(-2B(n)\right)\right)$.
We have $\frac{2}{5}\epsilon(n)^{1/3}\sim O(n^{-\xi})$ and $\frac{2}{5}e^{2}\cdot\exp\left(-2B(n)\right)\sim O(n^{-\xi})$,
so that $m(f,n)\sim O(n^{-\xi})$. Furthermore, $\delta(n)\sim O(n^{-\xi/2})\rightarrow0$.
The assumptions on the density function in Theorem \ref{Noncompact Theorem}
are therefore satisfied.
\end{example}

The essence of the examples above is embodied in the following general result whose proof follows along the path suggested by the examples.
\begin{proposition}
\label{prop:Polynomial-Exp tail}Assume that a continuous density
supported on $\mathbb{R}^{d}$ is of the form (\ref{eq:poly tail eq})
or (\ref{eq:exp tail eq}) with a compact zero-density region $S_{0}$.
For fixed $0<\eta<\frac{1}{d},\,0<\xi<\frac{1-2\eta d}{2}$, we can
find a sequence $\epsilon(n)\sim O(n^{-\psi}),\,0<\psi\leq\eta$, and
a decreasing sequence of $\delta(n)$, with $\lim_{n\rightarrow\infty}\delta(n)=0$,
such that 
\begin{align*}
M(f) & \coloneqq\sup_{\boldsymbol{w}\in\mathbb{R}^{d}}f(\boldsymbol{w})\in(0,\infty),\\
m(f,n)\coloneqq\min_{\boldsymbol{w}\in S_{f}^{1-\delta(n)}\backslash B_{\epsilon(n)}(S_{0})}f(\boldsymbol{w}) & \in(0,\infty)\sim O(n^{-\xi}),\,0<\xi<\frac{1-2\eta d}{2}.
\end{align*}
\end{proposition}

\begin{proof}
That $M(f)\coloneqq\sup_{\boldsymbol{w}\in\mathbb{R}^{d}}f(\boldsymbol{w})\in(0,\infty)$
is a direct consequence of the assumptions about the form of (\ref{eq:poly tail eq})
or (\ref{eq:exp tail eq}). 

We want to show that we can construct a decreasing sequence $\delta(n),$
with $\lim_{n\rightarrow\infty}\delta(n)=0$, such that the corresponding
$1-\delta(n)$ supports $S_{f}^{1-\delta(n)}\coloneqq[-B(n),B(n)]^{d}$
have minimal values $m(f,n)\sim O(n^{-\xi})$. First, we choose the
sequences $B(n)$ and $\epsilon(n)$, which jointly determine the desired
decay rate of $m(f,n)$. Then, if $B(n)$ is an increasing sequence,
by the definition of $S_{f}^{1-\delta(n)}$ as a cube, $\delta(n)$
is automatically a decreasing sequence. 

To determine $B(n)$, we use a sandwich argument on the minimal value
$m(f,n)$ attained on $S_{f}^{1-\delta(n)}$. We will consider two
sequences of balls, one where each ball is contained in $S_{f}^{1-\delta(n)}$,
called the sequence of \emph{ inner tangential balls} $\mathcal{B}_{*}(n)$,
and the other where each ball contains $S_{f}^{1-\delta(n)}$, called
the sequence of \emph{outer inclusive balls} $\mathcal{B}^{*}(n)$: 
\begin{align}
\mathcal{B}_{*}(n) & \coloneqq B_{B(n)}(\boldsymbol{0})\equiv\{\boldsymbol{x}\in\mathbb{R}^{d}\mid d(\boldsymbol{x},\boldsymbol{0})<B(n)\},\\
\mathcal{B}^{*}(n) & \coloneqq B_{4\sqrt{d}B(n)}(\boldsymbol{0})\equiv\{\boldsymbol{x}\in\mathbb{R}^{d}\mid d(\boldsymbol{x},\boldsymbol{0})<4\sqrt{d}B(n)\}.
\end{align}
Consider the following minimal values of $f$ on $\mathcal{B}_{*}(n)\backslash B_{\epsilon(n)}(S_{0})$
and $\mathcal{B}^{*}(n)\backslash B_{\epsilon(n)}(S_{0})$, 
\begin{align*}
m_{*}(f,n) & =\inf_{\boldsymbol{w}\in\mathcal{B}_{*}(n)\backslash B_{\epsilon(n)}(S_{0})}f(\boldsymbol{w}),\\
m^{*}(f,n) & =\inf_{\boldsymbol{w}\in\mathcal{B}^{*}(n)\backslash B_{\epsilon(n)}(S_{0})}f(\boldsymbol{w}).
\end{align*}
 Since $\mathcal{B}_{*}(n)\subset S_{f}^{1-\delta(n)}\subset\mathcal{B}^{*}(n)$
by the above construction, we have $m_{*}(f,n)\geq m(f,n)\geq m^{*}(f,n)$.
When $f$ is of the form (\ref{eq:poly tail eq}) or (\ref{eq:exp tail eq}),
we can choose the sequence $\epsilon(n)=\left(\frac{1}{C_{1}}\cdot n^{-\xi}\right)^{1/\gamma_{1}}$
where $\gamma_{1}\geq\gamma>0$ such that $0<\frac{\xi}{\gamma_{1}}\leq\eta$.
We use this choice of $\epsilon(n)$ to ensure that $\epsilon(n)\sim O(n^{-\psi})\sim O(n^{-\frac{\xi}{\gamma_{1}}})$,
and $0<\psi\leq\eta$ hold. Note that $\epsilon(n)$ is a decreasing
sequence by definition. 

\textbf{(Polynomial tail)} When $f$ is of the form (\ref{eq:poly tail eq}),
we use the $\epsilon(n)$ above and choose another sequence $B(n)=\left(\frac{1}{C_{2}}\cdot n^{\xi}\right)^{-1/\chi}$
and take $n$ sufficiently large so that $\epsilon(n)<\frac{1}{4}B(n)$
and $\epsilon_{0}<B(n)$. Then the minimal values for $\mathcal{B}_{*}(n)$ satisfy 
\begin{flalign*}
m_{*}(f,n) & \leq\min\left(C_{1}\cdot\epsilon(n)^{\gamma},C_{2}\cdot(B(n)-\epsilon(n))^{\chi},C_{2}\cdot B(n)^{\chi}\right)\\
 & \sim O\left(\min\left(C_{1}\cdot\epsilon(n)^{\gamma},C_{2}\cdot B(n)^{\chi}\right)\right)\\
 & \sim O\left(\min\left(C_{1}\cdot n^{-\frac{\gamma}{\gamma_{1}}\xi},n^{-\xi}\right)\right).
\end{flalign*}
Since $\gamma/\gamma_{1}\leq1$  for $n>1$ large enough,
we have that this asymptotic order of magnitude is determined by the behavior of the second term and
$m_{*}(f,n)\sim O(n^{-\xi})$.
Similarly, the minimal values for $\mathcal{B}^{*}(n)$ satisfy 
\begin{flalign*}
m^{*}(f,n) & \geq\min\left(C_{1}\cdot\epsilon(n)^{\gamma},C_{2}\cdot(4\sqrt{d}B(n)+\epsilon(n))^{\chi}\right)\\
 & \sim O\left(\min\left(C_{1}\cdot\epsilon(n)^{\gamma},C_{2}\cdot B(n)^{\chi}\right)\right) 
 \\
 & \sim O(n^{-\xi}).
\end{flalign*}

\textbf{(Exponential tail)} When $f$ is of the form (\ref{eq:exp tail eq}),
we use the $\epsilon(n)$ above and choose another sequence $B(n)=\frac{1}{\beta}\log\left(\frac{1}{C_{2}}\cdot n^{-\xi}\right)$
and take $n$ sufficiently large so that $\epsilon(n)<\frac{1}{4}B(n)$
and $\epsilon_{0}<B(n)$. Then the minimal values for $\mathcal{B}_{*}(n)$ satisfy
\begin{flalign*}
m_{*}(f,n) & \leq\min\left(C_{1}\cdot\epsilon(n)^{\gamma},C_{2}\cdot\exp(\beta(B(n)-\epsilon(n)))\right)\\
 & \sim O\left(\min\left(C_{1}\cdot\epsilon(n)^{\gamma},C_{2}\cdot\exp(\beta B(n))\right)\right)\\
 & \sim O\left(\min\left(C_{1}\cdot n^{-\frac{\gamma}{\gamma_{1}}\xi},n^{-\xi}\right)\right).
\end{flalign*}
Since $\gamma/\gamma_{1}\leq1$  for $n>1$ large enough,
we have that this asymptotic order of magnitude is determined by the behavior of the second term and
$m_{*}(f,n)\sim O(n^{-\xi})$.
Similarly, the minimal values for $\mathcal{B}^{*}(n)$ satisfy 
\begin{flalign*}
m^{*}(f,n) & \geq\min\left(C_{1}\cdot\epsilon(n)^{\gamma},C_{2}\cdot\exp(\beta(4\sqrt{d}B(n)+\epsilon(n)))\right)\\
 & \sim O\left(\min\left(C_{1}\cdot\epsilon(n)^{\gamma},C_{2}\cdot\exp\left(\beta B(n)\right)\right)\right) 
 \\
 & \sim O(n^{-\xi}).
\end{flalign*}

By a sandwich argument with $\mathcal{B}_{*}(n)\subset S_{f}^{1-\delta(n)}\subset\mathcal{B}^{*}(n)$
and $m_{*}(f,n)\geq m(f,n)\geq m^{*}(f,n)$, we know that we can find
some $\frac{1}{4\sqrt{d}}\cdot\left(\frac{1}{C_{2}}\cdot n^{\xi}\right)^{-1/\chi}\leq B(n)\leq\left(\frac{1}{C_{2}}\cdot n^{\xi}\right)^{-1/\chi}$
for a density of the form (\ref{eq:poly tail eq}) or some $\frac{1}{4\sqrt{d}}\cdot\frac{1}{\beta}\log\left(\frac{1}{C_{2}}\cdot n^{-\xi}\right)\leq B(n)\leq\frac{1}{\beta}\log\left(\frac{1}{C_{2}}\cdot n^{-\xi}\right)$
for a density of the form (\ref{eq:exp tail eq}). This choice of an
increasing sequence of  $B(n)$ proves our claim that we can find a
sequence of $\delta(n),\,\text{with }\lim_{n\rightarrow\infty}\delta(n)=0$,
such that $M(f)\in(0,\infty)$ and $m(f,n) \sim O(n^{-\xi})$.
\end{proof}
\begin{corollary}
\label{cor:tails}
Consider a continuous density supported on $\mathbb{R}^{d}$ of the
form (\ref{eq:poly tail eq}) or (\ref{eq:exp tail eq}). We can find
a sequence $\epsilon(n)\sim O(n^{-\psi}),\,0<\psi\leq\eta$, and a
decreasing sequence of $\delta(n)$ values, with $\lim_{n\rightarrow\infty}\delta(n)=0$,
such that the associated sequence of $(1-\delta(n))$-supports $S_{f}^{1-\delta(n)}$
have corresponding coverings $\mathcal{B}_{r(n)}^{d}$ satisfying
A.5' and A.6'.  Then

\noindent
(A) If $\eta$ and $\psi$ satisfy 
\begin{alignat*}{1}
\begin{array}{c}
1-2\eta d-2\overline{K_{f}}\cdot\psi>0,\end{array}
\end{alignat*}
we have $\lim_{n\rightarrow\infty}P(\text{\hluo{there are no empty }\ensuremath{\epsilon}(n)-outside balls})=1$.

\noindent
(B) If $\eta$ and $\psi$ satisfy 
\begin{align*}
\begin{array}{c}
1+d_{0}\eta-\underline{K_{f}}\eta-d\eta<0,\end{array}
\end{align*}
we have $\lim_{n\rightarrow\infty}P(\text{all \ensuremath{\epsilon}(n)-inside balls are empty})=1$. 
\end{corollary}

\begin{proof}
By the result in Proposition \ref{prop:Polynomial-Exp tail}, it suffices
to consider the grid construction in Lemma~\ref{lem:Exist regular covering}.
The cardinality of the resulting coverings $\mathcal{B}_{r(n)}^{d}$ of $S_{f}^{1-\delta(n)}$
can be guaranteed to be $O\left(\frac{B(n)^{d}}{r(n)^{d}}\right)$.  So, for case~(\ref{eq:poly tail eq}) we can have $\left|\mathcal{B}_{r(n)}^{d}\right|\sim$ 
\newline$O\left(\left(\frac{1}{C_{2}}\cdot n^{-\xi}\right)^{d/\chi}\cdot n^{\eta d}\right)$ and for case~(\ref{eq:exp tail eq}) we can have $\left|\mathcal{B}_{r(n)}^{d}\right|\sim$
\newline$O\left(\frac{1}{\beta^{d}}\log^{d}\left(\frac{1}{C_{2}}\cdot n^{-\xi}\right)\cdot n^{\eta d}\right)$.
In either case, $\left|\mathcal{B}_{r(n)}^{d}\right|$ is obviously
$o(n^{\Omega})$ for some sufficiently large $\Omega>0$. Therefore
all assumptions in Theorems \ref{Decreasing sequence of deltas} and
\ref{Decreasing sequence of deltas-1} are met.
\end{proof}

\section{\label{sec:Simulations-and-Connections}Simulations and Connections
to Other Areas}

\subsection{Simulation Studies for the choice of \texorpdfstring{$r,\epsilon$}{r,epsilon}}

Our theoretical results establish the rate of decay for the radius $r(n)$ and the neighborhood width $\epsilon(n)$. That is, for any positive constants $M_{r}$ and $M_{\epsilon}$,
\begin{align*}
r(n) & =M_{r}\cdot O(n^{-\eta}),\\
\epsilon(n) & =M_{\epsilon}\cdot O(n^{-\psi}).
\end{align*}
(with $\eta$ and $\psi$ following the same notation used in Theorem \ref{Main Theorem}) the asymptotic guarantees of filled and empty balls hold.  These guarantees allow one to identify lower dimensional zero-density regions in the limit.  However, for a fixed sample size, the actual values of $r(n)$ and $\epsilon(n)$ matter.  This subsection reports simulations investigating choice of the constants $M_r$ and $M_{\epsilon}$.

We consider the density
$f(\boldsymbol{x})\propto d(\boldsymbol{x},S_{0})^{4}\circ\boldsymbol{1}_{[0,1]^{2}}$.  
The zero-density region $S_0 = \{\frac{1}{2} \} \times [0.25,0.75]$ is of strictly lower dimension. For this example, $m(f,n)$ decays as a polynomial with $\xi=4\cdot\psi$ and $M(f) < 2$.  The conditions of Theorem \ref{Main Theorem}, parts (A) and (B), hold with, for example, $\eta=0.21$ and $\psi=0.01$ (therefore $0<\xi\leq\frac{1-4\eta}{2}$).

In the simulation, random samples of various sizes are generated from the density.  A grid of balls with centers on a lattice is used to cover the unit square.  The centers of the balls depend on $M_r$ and $n$ and follow the rule we described in Proposition \ref{prop:correct-covering-for-line-segment}. We note that the fraction of filled balls of a particular kind serves as an empirical estimate of the mean probability that a particular kind of ball is non-empty.  

\begin{figure}[ht!]
\adjustbox{center}{
\includegraphics[scale=0.4]{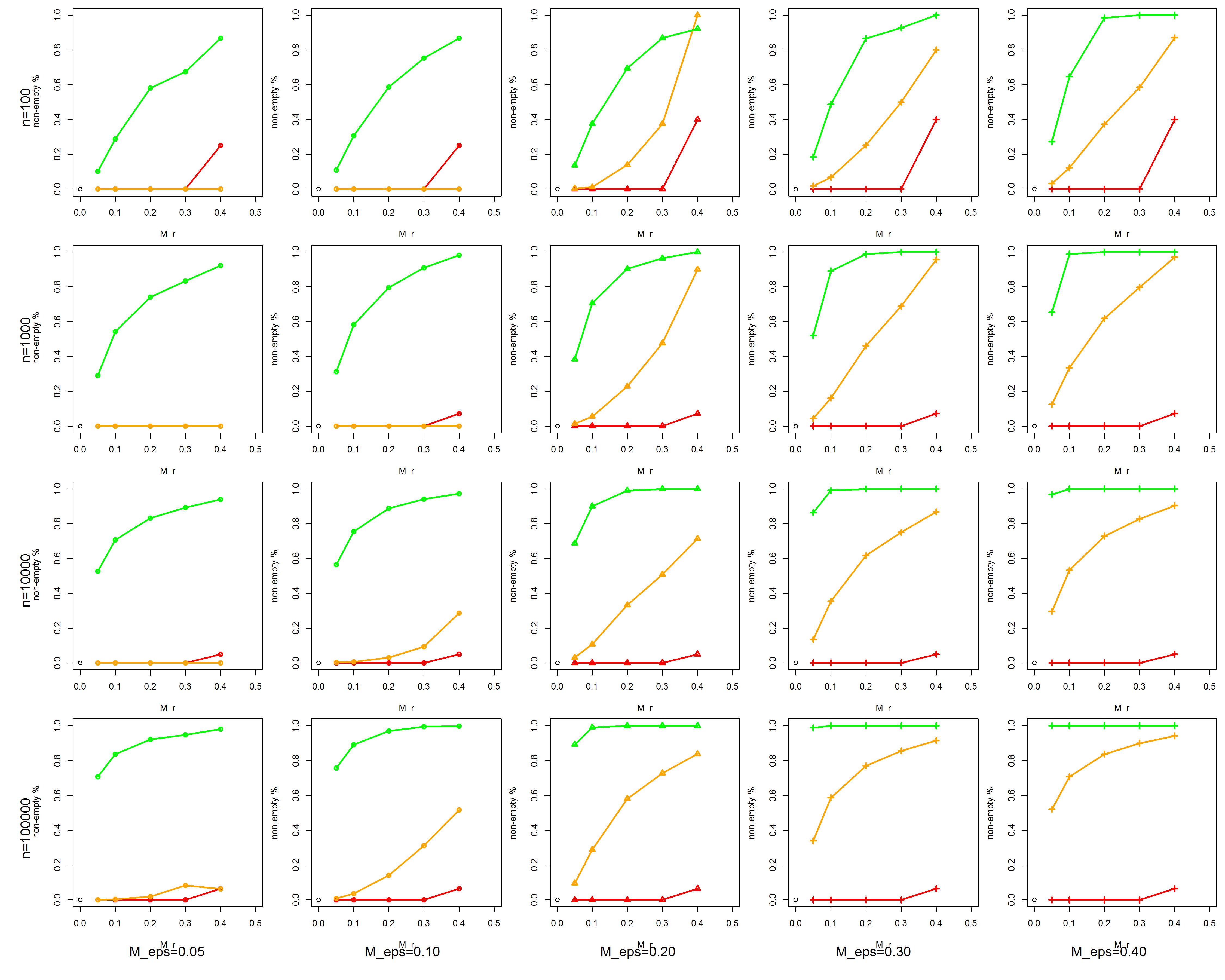}
}
\caption{\label{Summary Plot}The percentage of various types of non-empty covering balls
under different multipliers $M_{r}$ and $M_{\epsilon}$ for ball radius
$r(n)$ of the covering and neighborhood width of $S_0$, $\epsilon(n)$, for various sample sizes $n$.\protect \\ 
}
\end{figure}

Figure~\ref{Summary Plot} presents the results of the simulation. The color scheme is the same as in Figure \ref{Different kinds of balls}.
Proceeding down a column, the sample size changes from 100 to 10,000. Tracking the red lines, we see that the percentage of filled $\epsilon$-inside balls tends to $0$ as $n$ increases.  The orange lines represent the $\epsilon$-neighboring balls. The percentage of these that are filled does not always tend to $0$.  The green lines represent the $\epsilon$-outside balls.  Here, the fraction filled increases toward the eventual limit of~$1$.

As we step across rows of the figure, the value of $M_{\epsilon}$ increases from $0.05$ to $0.40$.  Within a plot, the value of $M_r$ increases from $0.05$ to $0.40$.  As $r(n)$ increases, the balls are larger and a greater fraction are filled, as evidenced by the lines within each plot.  As $M_{\epsilon}$ increases, balls in the covering are moved from $\epsilon$-outside to $\epsilon$-neighboring to $\epsilon$-inside.  The changing values of $\epsilon(n)$ and $r(n)$ shuffle the group membership of the balls and change the resulting fractions.  In general, the larger values of $M_r$ that we have investigated lead to greater separation of the red and green lines--hence greater differentiation between $\epsilon$-inside and $\epsilon$-outside balls.  Larger values of $M_{\epsilon}$ have a similar, though weaker effect.

The simulation also provides a caution.  The plots toward the bottom of the figure show much greater separation between the red and green lines.  This separation (near 0, near 1) is needed in order to reliably detect a zero-density region. For this density, the larger sample sizes prove much more effective than do the smaller sample sizes. 

The simulation results suggest an interesting possibility---the use of multiple coverings with balls of different sizes.  Doing so could allow one to examine the set of $M_r$, $M_{\epsilon}$ pairs that suggest the presence of a zero-density region. In practice, when $n$ is finite, the choice of multipliers $M_r$ and $M_{\epsilon}$ impacts the results.

\subsection{\label{KDE section}Connections to Existing Research\label{ML section}}

Our object of interest in this paper is the inverse image  of $\{0\}$ under the continuous density function
$f$\steve{, when it is of lower dimension than the support of $f$.}

\hluo{
The lower dimensional zero-density region is of interest to the topological data analysis (TDA) research community. Topological data analysis methods have been developed to explore and understand the topological structure of  the space where the data arise (for example, the support of a probability distribution) using a finite amount of observed data, as introduced in \cite{Carlsson2009} and \cite{Edelsbrunner&Harer2010}.

Having the data in hand, it is quite natural to consider them as a
random sample from some probability distribution. Interest then focuses
on the limiting behavior of geometric complexes based on a point set of growing size. In this regard, \cite{kahle2009topology} and \cite{Kahle&Meckes2013}
have developed results about the limiting behavior of Betti numbers
of complexes, based on the probabilistic results for random geometric
graphs given by \cite{Penrose2003}. Subsequent work in this direction
characterizes other types of topological features \cite{Duy2016,Bauer&Pausinger2018,kalisnik2019geometric}.
These results provide a picture of how the probability mechanism informs
us about the underlying topology. \cite{Bobrowski&Kahle2017} provide
a comprehensive review of the results along this line. This body of literature suggests that asymptotic regimes provide a convenient framework for analyzing the topology of data,}
\mario{as we do in the current paper.}

\hluo{Another line of relevant literature is on the set estimation}, \cite{Cuevas_2009}, \cite{Rigollet&Vert2009},  and \cite{singh2009adaptive} study the related concept of the  \emph{$\lambda$-level set for density $f$}. That is, the
set $S_{\lambda}\coloneqq\left\{ \boldsymbol{x}\in\mathbb{R}^{d}\mid f(\boldsymbol{x})\geq\lambda\right\},~\lambda>0$.
The special case where $\lambda=0$ is studied in \cite{devroye1980detection}. The difference between the regions is that $S_{\lambda}$  is the inverse image of an open set  that contains positive mass if non-empty while $S_0\coloneqq\left\{ \boldsymbol{x}\in\mathbb{R}^{d}\mid f(\boldsymbol{x})=0\right\} $ is the inverse image of a closed set and contains no mass. To the best of our knowledge, regions like $S_{0}$ have not previously been studied.

In the level set estimation literature, there are two major approaches to the problem of estimating level sets. Suppose we are interested in estimating the level sets of density $f$ (corresponding to a probability measure $F$). One approach is to construct
plug-in estimators $\hat{S_{\lambda}}\coloneqq \left\{ \boldsymbol{x}\in\mathbb{R}^{d}\mid\hat{f}(\boldsymbol{x})>\lambda\right\} $
for a level set, based on a kernel density estimator $\hat{f}$ computed
from the data and appropriate choices of the bandwidth parameters for the kernel \cite{Rigollet&Vert2009,Mason&Polonik2009, rinaldo2010generalized}. 
\hrluonew{This approach connects to the topology of the density when the level, $\lambda$, is varied \citep{chazal2013persistence,Fasy2014b}.  Asymptotically, a consistent estimator of the density may reveal persistent features.}

The other approach is based on the empirical excess-mass functional. The (empirical) excess-mass
functional $E(\lambda)\coloneqq F(S_\lambda)-\lambda \text{Leb}(S_\lambda)$ measures how the ``excess probability mass'' of the probability measure $F$ distributes over the region $S_\lambda$ when compared to Lebesgue measure \cite{hartigan1987estimation}. If we substitute the set $S_\lambda$ with a set estimator $\hat{S_{\lambda}}$ for the $\lambda$-level set $S_{\lambda}$, we can similarly consider the functional  $H(\lambda)\coloneqq F(\hat{S_{\lambda}})-\lambda \text{Leb}(\hat{S_{\lambda}})$, which can be used to evaluates the estimator $\hat{S_{\lambda}}$. The functional $H$ is maximized over a class of sets to obtain the level set estimator $\hat{S_{\lambda}}$ \cite{Polonik1995,walther1997granulometric}, to obtain level set estimators $\hat{S_{\lambda}}$. The
consistency and asymptotic behavior of both approaches has been derived
under regularity assumptions. 

In this paper, we study the lower dimensional
object that arises as the inverse image of \hluo{the} density function $f$ of a single point set. When we restrict our concern to the manifold $supp(f)$ and the
inverse of the density function $f^{-1}$ is sufficiently smooth,
$f^{-1}(\{a\}),\,f^{-1}(\{b\})$ encode the boundary of the manifold defined
by the inverse image $f^{-1}((a,b))$, as a sub-manifold of $supp(f)$. When $a=\lambda$ and $b=\infty$,
this inverse image defines $\lambda$-level sets. The object $S_{0}=f^{-1}(\{0\})$
is the boundary of such a particular example where $a=0$ and $b=\infty$. 
Unlike the topological features associated with the level sets, the feature we study does not necessarily generate equivalent classes in (persistent) homology. However, in applications like image segmentation, the boundaries are features of major interest \cite{LJS2019}.

Our method detects this specific kind of lower dimensional topological feature that could arise as a manifold boundary. We construct a sequence
of covering families to detect $S_0$ and we derive a set of sufficient
conditions that ensure consistent detection \hluo{in an asymptotic sense}. 
As is to be expected, when more sample points are available, our covering scheme locates the zero-density region
$S_{0}$ more accurately. In applied scenarios where the boundaries arise as zero-density regions of certain density functions, our method could help in detection \cite{LJS2019}.

\hrluonew{We do not provide an estimator of $S_0$ in this paper. A referee suggested that the union of empty balls in a covering could constitute an estimator of $S_0$. Additionally, if $M$ were not known, the same referee suggested that our method could be applied to a preliminary estimate of $M$. This is an interesting direction, but it is too difficult to develop fully in this article. For one thing, the performance of the method would depend on the accuracy of the preliminary estimator. In addition, the theoretical properties of the method would become more intricate, owing to the dependence between the preliminary estimation of $M$ and the covering procedure.} 

The main results in this paper also exhibit the relation between topological features that arise as $S_{0}$, and its dimensionality. In \cite{Adler_etal2013} and \cite{Owada&Adler2017}, the authors observe that higher dimensional topological features are generally smaller in scale. As shown by the sufficient conditions in Theorem \ref{Main Theorem} and \ref{Noncompact Theorem},
when the dimension $d_{0}$ of $S_{0}$ is higher, we need to specify
a faster decay rate for radii $r(n)$ of the covering sequence in
order to meet the sufficient conditions. This interplay between the dimension of the support, the dimension of the zero-density region $S_0$ and the sufficient growth rates supports \cite{Adler_etal2013}'s observation from a different angle. 

\section{\label{sec:Conclusion}Conclusions}

In this paper, we consider the question of detection of $S_{0}=f^{-1}(\{0\})$ lying in the support $M$ of a continuous density $f$ when $S_0$ has Minkowski dimension strictly lower than the dimension of the data, $d$. This type of topological feature is difficult to identify and has not been studied before. The main contribution in this paper is to provide a novel approach, based on a sequence of coverings tied to growing sample sizes, to study a specific kind of lower dimensional object $S_0$ in the support of a density function. This approach works under both compact and non-compact support and its construction has geometric intuition. Being of lower dimension, $S_0$ is a delicate object. It is in the closure of $M\backslash S_0$ and, in a sense, disappears when one looks at it at too coarse a resolution. 

Our strategy is to construct a sequence of covering balls (e.g., Lemma
\ref{lem:Exist regular covering}) and shrink their radius as the
sample size $n$ goes to infinity. We derive a set of sufficient conditions,
using the local behavior of $f$, captured by 
$\overline{K_{f}}$ and $\underline{K_{f}}$,
that  ensure that a particular covering scheme leads to probability one consistency results in
Theorem \ref{Main Theorem} (compact support) and Theorem \ref{Noncompact Theorem}
(non-compact support). This set of theorems in Section \ref{sec:Statement-of-the-main-result}
and \ref{sec:Non-compact-support} can be generalized to the case where
$S_{0}$ has multiple  disconnected components.  
As the sample size $n$ tends to $\infty$,  a shrinking  $\epsilon$-neighborhood
of $S_{0}$ can be identified by empty covering balls asymptotically
with probability one while the rest of $M$ will be covered by non-empty
balls. 

Our result provides a range of asymptotic schemes for the radius of
covering balls under which the lower dimensional topological feature
can be detected. We provide experimental evidence to support
our claim in Section \ref{sec:Simulations-and-Connections}.
Our approach supports the insight that different dimensional
topological features occur at different scales. 
The novelty of our result is the role of the ambient dimension of the data in addition to the local behavior of $f$ near $S_0$.

Our approach and results focus on the connection between i.i.d.\ samples and the near-topological features of the support of the density from which they are drawn.  As future work, it is of great interest to generalize the results to dependent draws from the density, with the eventual goal of understanding how probabilistic dependence in the sample can be useful in the construction of complexes based on sub-samples (e.g. witness complexes) used in TDA \cite{HLuoGeneralized2020,hrluo_2021d}.

\begin{appendix}
\section{Proof of Main Theorem}\label{Proof of the Main Theorem}

Theorem \ref{Main Theorem} is established by providing a bound on the probability
mass of a covering ball $B_{r}(\boldsymbol{x})$, counting the number
of each type of covering ball under consideration and taking a limit
as the sample size $n\rightarrow\infty$. A sequence of upper bounds is needed for
$\epsilon$-inside balls and a sequence of lower bounds is needed
for $\epsilon$-outside balls.

In the statement of the theorem, we assume that the sequence 
\begin{align*}
r(n)\sim O(n^{-\eta}) & ,\,0<\eta<\frac{1}{d},\\
\epsilon(n)\sim O(n^{-\psi}) & ,\,\psi\leq\eta,\\
\text{and }2r & (n)\leq\epsilon(n)<1.
\end{align*}
To establish 
\steve{a bound} on the probability mass
$P(B_{r}(\boldsymbol{x}))=\int_{B_{r}(\boldsymbol{x})}f(\boldsymbol{w})d\boldsymbol{w}$
of a ball $B_{r}(\boldsymbol{x})$, it suffices to consider the volume
of the ball and a bound on the density function over the ball. The
volume for a $d$-dimensional ball of radius $r$ is 
\begin{align}
V_{d}(r) & ={\displaystyle \frac{\pi^{\frac{d}{2}}}{\Gamma\left(\frac{d}{2}+1\right)}r^{d},}\label{eq:volume of ball}
\end{align}
which for our sequence of radii $r(n)$ is $V_{d}(r(n))\sim O(n^{-d\eta})$.
Recall from the definitions of $\epsilon$-inside and $\epsilon$-outside
covering balls:

\emph{$\epsilon$-inside balls} are those balls $B_{r}(\boldsymbol{x})$
such that $\boldsymbol{x}\in B_{\epsilon}(S_{0})$ and $B_{r}(\boldsymbol{x})\cap S_{0}\neq\emptyset$.

\emph{$\epsilon$-outside balls} are those balls $B_{r}(\boldsymbol{x})$
such that $\boldsymbol{x}\notin B_{\epsilon}(S_{0})$.  

Outside $B_{\epsilon(n)}(S_{0})$ the density is bounded below by
$m(f,n)\sim O(n^{-\xi})$. Inside $B_{\epsilon(n)}(S_{0})$ the density
is bounded below by the inequality in the \hyperlink{Assumption3}{Assumption
A.4}:
\begin{align}
L_{f}\cdot d(\boldsymbol{x},S_{0})^{\overline{K_{f}}(\epsilon(n))} & \leq f(\boldsymbol{x})\leq U_{f}\cdot d(\boldsymbol{x},S_{0}){}^{\underline{K_{f}}(\epsilon(n))}.\label{eq:ord-control}
\end{align}

\begin{itemize}
\item \emph{Upper bound for an $\epsilon(n)$-inside covering ball.}\\
 Since the $\epsilon(n)$-inside ball intersects $S_{0}$, i.e. $B_{r(n)}(\boldsymbol{x})\cap S_{0}\neq\emptyset$,
we know that any point $\boldsymbol{y}\in B_{r(n)}(\boldsymbol{x})$
will be at most $2r(n)\leq\epsilon(n)$ away from the zero-density
region $S_{0}$. The probability mass contained in an $\epsilon(n)$-inside
covering ball $B_{r(n)}(\boldsymbol{x})$ is bounded from above by
\begin{align}
P(B_{r(n)}(\boldsymbol{x}))\leq & U_{f}\cdot V_{d}(r(n))\cdot(2r(n))^{\underline{K_{f}}(\epsilon(n))}
\end{align} 
There are two types of $\epsilon(n)$-outside balls: those that are
entirely contained in $M=[0,1]^{d}$ and those that are only partially
contained in $M$. Probability mass $P(B_{r(n)}(\boldsymbol{x}))=\int_{B_{r(n)}(\boldsymbol{x})}f(\boldsymbol{w})d\boldsymbol{w}$
is bounded from below by \steve{the} volume of the ball times the minimum of the
density $f(\boldsymbol{x})$ in the ball. 
\item \emph{Lower bound for an $\epsilon(n)$-outside covering ball that lies
within $M$. }\\
 From the assumption on density $f$,
\begin{align*}
m(f,n)\coloneqq\min_{\boldsymbol{w}\in M\backslash B_{\epsilon}(S_{0})}f(\boldsymbol{w}) & \in(0,\infty)\sim O(n^{-\xi}),\,0<\xi<\frac{1-2\eta d}{2}.
\end{align*}
The probability mass contained in the $\epsilon(n)$-outside covering
ball is bounded from below by 
\begin{align}
P(B_{r(n)}(\boldsymbol{x}))\geq V_{d}(r(n))\cdot\min\left[L_{f}\cdot(\epsilon(n)-r(n))^{\overline{K_{f}}(\epsilon(n))},m(f,n)\right] & .\label{eq:lower bound2}
\end{align}
\item \emph{Lower bound for an $\epsilon(n)$-outside covering ball that does
not lie entirely within $M$. }\\
 \hyperlink{Assumption3}{Assumption A.5} states that the center
of the covering ball $B_{r(n)}(\boldsymbol{x})$ is in $M$. We know
that the volume of such an $\epsilon(n)$-outside ball $B_{r(n)}(\boldsymbol{x})$
will be at least $(\frac{1}{2})^{d}$ times the volume $V_{d}(r(n))$
of a full $d$-dimensional ball. With the same $m(f,n)>0$ we have
\begin{align}
P(B_{r(n)}(\boldsymbol{x}))\geq\left(\frac{1}{2}\right)^{d}V_{d}(r(n))\cdot\min\left[L_{f}\cdot(\epsilon(n)-r(n))^{\overline{K_{f}}(\epsilon(n))},m(f,n)\right]\label{eq:lower bound3}
\end{align}
This lower bound also holds for $\epsilon(n)$-outside balls that
are entirely contained in $M$. 
\end{itemize}
In the following discussion, we use the notation $p_{B_{r(n)}(\boldsymbol{x})}\coloneqq P(B_{r(n)}(\boldsymbol{x}))$
to emphasize that we use it as a parameter.

\subsection*{Part (A) (Outside balls)}

Consider regular families of covering balls $\mathcal{B}_{r(n)}^{d}$
of $supp(f)$ and a sequence of $\epsilon(n)$-neighborhoods of the
zero-density region $S_{0}$ where $\dim S_{0}<d$. The Hoeffding
concentration inequality for $X\sim Binom(n,p)$ can be written for
arbitrary $\gamma>0$, 
\begin{align}
{\displaystyle P\left\{ (p-\gamma)n\leq X\leq(p+\gamma)n\right\} \geq1-2\exp\left(-2\gamma^{2}n\right).}\label{eq:Hoeffding}
\end{align}
The inequality ensures that the number of observations falling into
a single ball $B_{r(n)}(\boldsymbol{x})$ will be close to its expectation.
In a single ball, the number $N_{B_{r(n)}(\boldsymbol{x})}$ of observations
falling in the $d$-dimensional ball of radius $r(n)$, is distributed
as a binomial distribution $Binom(n,p_{B_{r(n)}(\boldsymbol{x})})$.
Also, 
\begin{align*}
P(B_{r(n)}(\boldsymbol{x})\text{ not empty}) & =P\left(N_{B_{r(n)}(\boldsymbol{x})}\geq1\right).
\end{align*}

By assumption on the range of $\eta$, the lower bounds (\ref{eq:lower bound2}),
(\ref{eq:lower bound3}) under the different situations are determined
by the factor $r(n)^{d}$. But note that $r(n)^{d}\sim O(n^{-d\eta})$
with $-d\eta>-1$ as assumed. We can choose $\mathsf{N}(\boldsymbol{x})\in\mathbb{N}$
large enough such that $\mathsf{N}(\boldsymbol{x})\cdot p_{B_{r(n)}(\boldsymbol{x})}\geq1$,
ensuring $n\cdot p_{B_{r(n)}(\boldsymbol{x})}\geq1$ when $n\geq\mathsf{N}(\boldsymbol{x})$.
We take $\mathsf{N}(\boldsymbol{x})=\frac{1}{p_{B_{r(n)}(\boldsymbol{x})}}$
and observe that $\mathsf{N}(\boldsymbol{x})\sim O(n^{d\eta})$ with
$d\eta<1$, rendering such a choice of $n\geq\mathsf{N}(\boldsymbol{x})$
possible. Assume in the arguments below that $n\geq\mathsf{N}(\boldsymbol{x})$.

For an $\epsilon(n)$-outside ball $B_{r(n)}(\boldsymbol{x})$, we
can rewrite the probability using (\ref{eq:Hoeffding}) with $X=N_{B_{r(n)}(\boldsymbol{x})},p=p_{B_{r(n)}(\boldsymbol{x})}$
and $\gamma=\frac{p_{B_{r(n)}(\boldsymbol{x})}}{2}$, 
\begin{align}
P(B_{r(n)}(\boldsymbol{x})\text{ not empty}) & \geq P\left(\left|N_{B_{r(n)}(\boldsymbol{x})}-np_{B_{r(n)}(\boldsymbol{x})}\right|<n\cdot\frac{p_{B_{r(n)}(\boldsymbol{x})}}{2}\right)\nonumber \\
 & \geq\left[1-2\exp\left(-\frac{1}{2}p_{B_{r(n)}(\boldsymbol{x})}^{2}\cdot n\right)\right].\label{eq:estimate}
\end{align}
From the discussion of the probability mass contained in each type
of covering ball above, if $B_{r(n)}(\boldsymbol{x})$ is an $\epsilon(n)$-outside
covering ball, then from formulas (\ref{eq:lower bound2}) and (\ref{eq:lower bound3})
we have 
\begin{align}
& P(B_{r(n)}(\boldsymbol{x})\text{ not empty}) \nonumber\\ & \geq\left[1-2\exp\left(-\frac{1}{2}\left[V_{d}(r(n))\cdot\min\left[L_{f}\cdot(\epsilon(n)-r(n))^{\overline{K_{f}}(\epsilon(n))},m(f,n)\right]\right]^{2}\cdot n\right)\right]\label{epsilon-outside ball bound-1}
\end{align}

If $\min\left[L_{f}\cdot(\epsilon(n)-r(n))^{\overline{K_{f}}(\epsilon(n))},m(f,n)\right]=m(f,n)$
then 
\begin{align*}
P(B_{r(n)}(\boldsymbol{x})\text{ not empty}) & \geq\left[1-2\exp\left(-\frac{1}{2}\left[V_{d}(r(n))\cdot m(f,n)\right]^{2}\cdot n\right)\right]
\end{align*}

If $\min\left[L_{f}\cdot(\epsilon(n)-r(n))^{\overline{K_{f}}(\epsilon(n))},m(f,n)\right]=L_{f}\cdot(\epsilon(n)-r(n))^{\overline{K_{f}}(\epsilon(n))}$
then 
\begin{align*}
& P(B_{r(n)}(\boldsymbol{x})\text{ not empty}) \nonumber \\ & \geq\left[1-2\exp\left(-\frac{1}{2}\left[V_{d}(r(n))\cdot L_{f}\cdot(\epsilon(n)-r(n))^{\overline{K_{f}}(\epsilon(n))}\right]^{2}\cdot n\right)\right]
\end{align*}
Consider the right hand side of (\ref{epsilon-outside ball bound-1}), taking
$n\rightarrow\infty$. If it holds that 
\begin{align}
\begin{array}{c}
1-2\eta d-2\xi>0\end{array} \text{and }  1-2\eta d-2\overline{K_{f}}\cdot\psi>0\label{eq:Sufficient A}
\end{align}
then the limit $\lim_{n\rightarrow\infty}P(B_{r(n)}(\boldsymbol{x})\text{ not empty})$
is 1. For the first inequality, it is simply $\xi<\frac{1-2\eta d}{2}\leq\frac{1}{2}$ (the second equality holds iff $d=0$)
as we assumed in the statement of the theorem so it suffices to consider
the second inequality.

We want to let $n>\frac{1}{p_{B_{r(n)}(\boldsymbol{x})}}$ for all
$\boldsymbol{x}\in M\backslash B_{\epsilon}(S_{0})$. But $\frac{1}{p_{B_{r(n)}(\boldsymbol{x})}}$
reaches its maximum when $p_{B_{r(n)}(\boldsymbol{x})}$ attains its
minimum $m(n,f)\coloneqq\min_{\boldsymbol{w}\in M\backslash B_{\epsilon}(S_{0})}f(\boldsymbol{w})$,
i.e. $\frac{1}{p_{B_{r(n)}(\boldsymbol{x})}}\leq\frac{1}{V^{d}(r(n))\cdot m(n,f)}\sim O(n^{\eta d+\xi})$.
As long as $\eta d+\xi<1$, which is guaranteed by (\ref{eq:Sufficient A}),
we can ensure that $n$ is greater than the maximum value $\max_{\boldsymbol{w}\in M\backslash B_{\epsilon(n)}(S_{0})}\frac{1}{p_{B_{r(n)}(\boldsymbol{w})}}$.
Therefore, such an $n$ is greater than all such quantities $\frac{1}{p_{B_{r(n)}(\boldsymbol{x})}},\forall\boldsymbol{x}\in M\backslash B_{\epsilon(n)}(S_{0})$
and so $n\cdot p_{B_{r(n)}(\boldsymbol{x})}\geq1$ holds simultaneously
for all $\boldsymbol{x}\in M\backslash B_{\epsilon(n)}(S_{0})$. 
\begin{align}
P(\text{\hluo{there are no empty }\ensuremath{\epsilon(n)}-outside balls}) & =P(\cap_{B\in\{\epsilon(n)-\text{outside balls}\}}\{B\text{ not empty}\})\nonumber \\
 & =1-P(\cup_{B\in\{\epsilon(n)-\text{outside balls}\}}\{B\text{ empty}\})\nonumber \\
 & \geq1-\sum_{B\in\{\epsilon(n)-\text{outside balls}\}}P(B\text{ empty})\nonumber \\
 & \geq1-\left|\mathcal{B}_{r(n)}^{d}\right|\cdot\sup_{B\in\{\epsilon(n)-\text{outside balls}\}}P(B\text{ empty})
\end{align}
By the equation (\ref{eq:estimate}) and the subsequent bound (\ref{epsilon-outside ball bound-1})
we obtain that 
\begin{align}
& P(\text{\hluo{there are no empty }\ensuremath{\epsilon(n)}-outside balls}) \nonumber \\
& \geq1-|\mathcal{B}_{r(n)}^{d}|\cdot\left(1-\inf_{B\in\{\epsilon(n)-\text{outside balls}\}}P(B\text{ not empty})\right)\nonumber \\
 & \geq1-(L_{1}\cdot n^{d\eta})\cdot\left(L_{2}\cdot\exp\left(-n^{1-2\eta d-2\overline{K_{f}}\cdot\psi}\right)\right)\label{eq:key to non-compact}
\end{align}
for sufficiently large $n$.The constants $L_{1},L_{2}\in(0,\infty)$ exist by the definition
of the big O notation for sufficiently large $n$. The last expression
follows from $|\mathcal{B}_{r(n)}^{d}|\sim O(n^{d\eta})$ and (\ref{epsilon-outside ball bound-1}).
It is dominated by the second term in the product. Therefore, if $1-2\eta d-2\overline{K_{f}}\cdot\psi>0$
then 
\begin{align*}
\lim_{n\rightarrow\infty}P(\text{\hluo{there are no empty }\ensuremath{\epsilon(n)}-outside balls}) & =1.
\end{align*}

\subsection*{Part (B) (Inside balls)}

We observe that the event ``all $B_{\epsilon(n)}(S_{0})$-inside
balls are empty'' can be regarded as all observations falling into
the other two types of covering balls. Note that, as stated in the
theorem, we only ensure that $\epsilon(n)$-inside covering balls
are empty but not every $\epsilon(n)$-outside ball contains at least
one observations under the assumptions in part (B). We investigate the
upper bound on the probability mass in each of these covering balls.

Since the observations are i.i.d.$\,$we have, using the notation
that \newline$p_{\cup_{B\in\{\epsilon(n)-\text{inside balls\}}}B}\coloneqq P(\cup_{B\in\{\epsilon(n)-\text{inside balls\}}}B)$,
\begin{align}
P(\text{all \ensuremath{\epsilon(n)}-inside balls are empty}) & =\left(1-p_{\cup_{B\in\{\epsilon(n)-\text{inside balls\}}}B}\right)^{n}\nonumber \\
 & \text{since covering balls may overlap,}\nonumber \\
 & \geq\left(1-\sum_{B\in\{\epsilon(n)-\text{inside balls\}}}p_{B}\right)^{n}
\end{align}
from the upper bound above.

\hyperlink{Assumption3}{Assumption 3} asserts that the limit $\underline{K_{f}}=\lim_{\epsilon\rightarrow0^{+}}\underline{K_{f}}(\epsilon)$
exists and that we can choose a bound $K_{U}\geq\lim_{\epsilon\rightarrow0^{+}}\underline{K_{f}}(\epsilon)$
uniformly. By the bounds (\ref{eq:ord-control}) for $\epsilon(n)$-inside
covering balls, there is a positive constant $D_{f,d}\coloneqq U_{f}\cdot\frac{\pi^{\frac{d}{2}}}{\Gamma\left(\frac{d}{2}+1\right)}\cdot2^{K_{U}}\geq U_{f}\cdot\frac{\pi^{\frac{d}{2}}}{\Gamma\left(\frac{d}{2}+1\right)}\cdot2^{\underline{K_{f}}(\epsilon(n))}$
depending only on density $f$ and the dimension $d$. The factor
$\frac{\pi^{\frac{d}{2}}}{\Gamma\left(\frac{d}{2}+1\right)}$ comes
from the multiplier of the volume of a $d$-dimensional ball (\ref{eq:volume of ball}).
We denote by $\mathcal{B}_{r(n)}^{d_{0}}$ the sub-collection of covering
balls from $\mathcal{B}_{r(n)}^{d}$ that intersect $S_{0}$. 
\begin{align}
& P(\text{all \ensuremath{\epsilon(n)}-inside balls are empty}) \nonumber \\
& \geq\left(1-|\mathcal{B}_{r(n)}^{d_{0}}|\cdot\left[U_{f}\cdot V_{d}(r(n))\cdot(2r(n)){}^{\underline{K_{f}}(\epsilon(n))}\right]\right)^{n}\nonumber \\
 & \geq\left(1-H_{\varepsilon}(n)\cdot\left[D_{f,d}\cdot(r(n)){}^{\underline{K_{f}}(\epsilon(n))+d}\right]\right)^{n}\label{eq:control2}
\end{align}
Due to the \hyperlink{Assumption5}{Assumption A.6}, the collection
of covering balls $\mathcal{B}_{r(n)}^{d_{0}}$ that intersect $S_{0}$
satisfy $|\mathcal{B}_{r(n)}^{d_{0}}|\leq H_{\varepsilon}(n)\sim O(n^{d_{0}\eta(1+\varepsilon)})$
for any $\varepsilon>0$. On one hand, by the assumption for part
(B), we have 
\begin{align}
\begin{array}{c}
1+d_{0}\eta-\underline{K_{f}}\eta-d\eta\end{array} & <0
\end{align}
and so we can find $1>\varepsilon>0$ small enough such that, 
\begin{align}
\begin{array}{c}
1+d_{0}\eta(1+\varepsilon)-\underline{K_{f}}\eta-d\eta\end{array} & <0\label{eq:vareps control}
\end{align}
We need a simple but nontrivial lemma known as the Bernoulli inequality
to proceed. 
\begin{lemma}
\label{(Bernoulli-inequality)}(Bernoulli inequality) If $1+x>0$,
then $(1+x)^{n}\geq1+nx$ for $\forall n\in\mathbb{N}$. 
\end{lemma}
\begin{proof}
Let us prove the lemma by induction. For $n=0,1$ this trivially holds
as an equality. Assume that the inequality $(1+x)^{n}\geq1+nx$ holds
for every $n\leq k\in\mathbb{N}$ and proceed by induction for $n=k+1$.
\begin{align*}
(1+x)^{n} & =(1+x)\cdot(1+x)^{n-1}=(1+x)\cdot(1+x)^{k}\\
 & \text{since }n-1=k\leq k\text{ by the inductive hypothesis,}\\
 & \geq(1+x)\cdot(1+kx)\\
 & \text{due to the fact that }(1+x)>0\text{ does not flip the inequality,}\\
 & =1+kx+x+kx^{2}=1+(k+1)\cdot x+kx^{2}\\
 & \geq1+(k+1)\cdot x,
\end{align*}

since $kx^{2}\geq0$,which completes the induction. 
\end{proof}
Let $x=H_{\varepsilon}(n)\cdot\left[D_{f,d}\cdot(r(n)){}^{\underline{K_{f}}(\epsilon(n))+d}\right]\sim O(n^{d_{0}\eta(1+\varepsilon)})\cdot\left[D_{f,d}\cdot(r(n)){}^{\underline{K_{f}}(\epsilon(n))+d}\right]$,
which is of order $O(n^{d_{0}\eta(1+\varepsilon)-\underline{K_{f}}\cdot\eta-d\eta})$
as $n\rightarrow\infty$. We emphasize again that $\varepsilon>0$ is an auxiliary quantity
which we first define in \hyperlink{Assumption5}{Assumption 5}.
The $\epsilon=\epsilon(n)$ is the quantity in Theorem \ref{Main Theorem}
that defines the neighborhood size of $S_{0}$. These two are different
quantities. Under the assumption 
\begin{align}
\begin{array}{c}
1+d_{0}\eta(1+\varepsilon)-\underline{K_{f}}\eta-d\eta<0\end{array} & \Leftrightarrow d_{0}\eta(1+\varepsilon)-\underline{K_{f}}\eta-d\eta<-1,
\end{align}
we can take $n$ sufficiently large, say $n\geq\mathsf{N}_{1}$, to
ensure $d_{0}\eta(1+\varepsilon)-\underline{K_{f}}\eta-d\eta<-1$
and therefore $1+x>0$. Then we can apply the Bernoulli inequality
to the right hand side of (\ref{eq:control2}) and we have the following lines,
with $H_{\varepsilon}(n)\sim O(n^{d_{0}\eta(1+\varepsilon)})$ and
(\ref{eq:control2}), as $n\rightarrow\infty$:
\begin{align}
P(\text{all \ensuremath{\epsilon(n)}-inside balls are empty}) & \geq\left(1-(L_{3}\cdot n^{d_{0}\eta(1+\varepsilon)})\cdot\left[D_{f,d}\cdot(r(n)){}^{\underline{K_{f}}+d}\right]\right)^{n}\nonumber \\
 & \geq\left(1-n\cdot(L_{3}\cdot n^{d_{0}\eta(1+\varepsilon)})\cdot\left[D_{f,d}\cdot(r(n)){}^{\underline{K_{f}}+d}\right]\right)\nonumber \\
 & \geq\left(1-L_{f,d}\cdot(n^{1+d_{0}\eta(1+\varepsilon)-\underline{K_{f}}\eta-d\eta})\right)\label{eq:control3}
\end{align}
for sufficiently large $n$.The constant $L_{3}\in(0,\infty)$ exists by the definition of the
big O notation for sufficiently large $n$. If we take the limit $n\geq\mathsf{N}_{1},n\rightarrow\infty$
on both sides of the inequality we know that the quantity $\left(1-L_{f,d}\cdot(n^{1+d_{0}\eta(1+\varepsilon)-\underline{K_{f}}\eta-d\eta})\right)$
converges to 1 from (\ref{eq:vareps control}). therefore if $1+d_{0}\eta-\underline{K_{f}}\eta-d\eta<0$
then 
\begin{align*}
\lim_{n\rightarrow\infty}P(\text{all \ensuremath{\epsilon(n)}-inside balls are empty}) & =1.
\end{align*}



\end{appendix}

\section*{Acknowledgements}
This material is based upon work supported by the National Science
Foundation under Grants No. DMS-1613110 and DMS-2015552, and No. SES-1921523.

\hrluonew{
We thank the anonymous referee, whose comments greatly improve the article. We thank the AE for helpful comments and handling. 
}

\bibliographystyle{siamplain} 
\bibliography{Draft_THESIS_bib}

\end{document}